\RequirePackage{fix-cm}
\documentclass[twocolumn, smallextended]{svjour3}
\smartqed  %
\usepackage{url}
\usepackage{natbib}
\usepackage{graphicx}
\usepackage{subfig}
\usepackage{newtxtext}
\usepackage{mathtools}
\usepackage{amsfonts}
\usepackage{dsfont}
\usepackage{rotating}
\usepackage{booktabs}
\usepackage{multirow}
\usepackage{array}
\usepackage{float}

\journalname{Here is the journal name}
\begin{document}

\title{Reliability-based topology optimization under random-field material model
}

\author{Trung Pham \and Christopher Hoyle}

\institute{Trung Pham \and Christopher Hoyle \at
              School of Mechanical, Industrial \& Manufacturing Engineering \\
              Oregon State University \\
              Corvallis, OR 97331--6001, USA \\
              \email{chris.hoyle@oregonstate.edu}
}

\date{Received: xx/xx/20xx / Accepted: date}

\maketitle

\begin{abstract}
This paper presents an algorithm for reliability-based topology optimization of linear elastic continua under random-field material model. The modelling random field is discretized into a small number of random variables, and then the interested output is estimated by a stochastic response surface. A single-loop inverse-reliability algorithm is applied to reduce computational cost of reliability analysis. Two common benchmark problems in literature are used for demonstration purposes. Different values of target reliability and ranges of Young's modulus are considered to investigate their effects on resulting optimized topologies. Lastly, Monte Carlo simulation tests the proposed algorithm for correctness and accuracy.
\keywords{Continuum topology optimization \and Material uncertainties \and Structural reliability \and Karhunen$-–$Lo\`{e}ve expansion \and Stochastic response surface}
\end{abstract}

\section{Introduction}
\label{intro}
In structural design, geometrical and topological properties of the design have significant impacts on structural performance. Topology optimization (TO) has emerged as an effective method to optimize such properties, including size, shape, and connectivity of the design \citep{bendsoe_topology_2003}. TO can identify designs with improved performance while using the least amount of material. However, deterministic inputs have often been assumed in research on TO, while most observable phenomena contain an inherent amount of uncertainty; in the analysis and design of engineering systems the presence of uncertainty has always induced an acute impact. Ignoring uncertainty may produce sub-optimal designs that perform poorly under random inputs. Among different sources of uncertainty, the material property is inherently random in space and can be modeled by a random field. Therefore, this paper proposes an algorithm for reliability-based topology optimization under a random-field material model. Specifically, the Karhunen$-–$Lo\`{e}ve expansion discretizes the modeling random field into a small number of random variables, and then the output of interest is estimated by the stochastic response surface method. The uncertainty in material property makes the problem constraints probabilistic, which will  require a drastically different approach for the optimization problem. A single-loop algorithm called Sequential Optimization and Reliability Assessment (SORA) \citep{du_sequential_2004} coupled with the performance measure approach \citep{tu_new_1999} is used to reduce the computational cost of reliability analysis. The proposed algorithm is tested using two numerical examples, and finally verified by Monte Carlo simulation.

The layout of this paper is as follows. Section \ref{background} is a
literature review on the current state of the art providing background of uncertainty propagation, reliability analysis, and TO under uncertainty. Section \ref{tou} states the mathematical formulation of the deterministic and reliability-based topology optimization (RBTO) problem. This section also explains in detail our proposed solution for the RBTO problem. Numerical examples, including a simply supported beam and a L-shape beam, are done in section \ref{results} with discussion of the results. Finally, the paper is concluded by findings and future work.

\section{Background}
\label{background}
Though there are different probability interpretations \citep{prob_inter}, in the context of continuum mechanics the assumption of continuity over the entire domain makes it suitable to apply the so-called Kolmogorov's probability theory \citep{kolmogorov_2018} to the problem of ``design under uncertainty''. Using such theory coupled with an appropriate solution method, uncertainty in TO has been tackled mainly in the settings of robust optimization (RO) and reliability-based optimization (RBO). 

In RO the primary goal is to minimize the variability of outputs of interest due to uncertainty \citep{taguchi_1986}. This is usually archived by optimizing a weighted sum of mean and standard deviation of the objective function. This approach was implemented in a number of papers with various sources of uncertainty and solution methods: spatial variation of manufacturing error with Monte Carlo simulations \citep{schevenels_robust_2011}; random-field truss material with a multi-objective approach \citep{richardson_robust_2015}; random loading field and random material field with the level set method \citep{chen_level_2010}; material and geometric uncertainties with stochastic collocation methods and perturbation techniques \citep{lazarov_topology_2012,lazarov_topology_2012-1}; misplacement of material and imperfect geometry \citep{jansen_robust_2013,jansen_robust_2015}; Young's modulus of truss members with a perturbation method \citep{asadpoure_robust_2011}; random-field material properties with a polynomial chaos expansion \citep{tootkaboni_topology_2012}; geometric and  material properties uncertainties with a stochastic perturbation method for frame structures \citep{changizi_robust_2017}; random loading field with stochastic collocation methods \citep{zhao_robust_2015}. A unified framework for robust topology optimization (RTO) was proposed in \cite{richardson_unified_2016}, while \cite{zhao_robust_2014} solved the RTO problem exploiting the linear elasticity of structure. Some authors \citep{jalalpour_efficient_2016,richardson_robust_2015} claimed that one weakness of the RO methodology is the seemingly arbitrary weighting factors, which actually reflect risk-taking attitude of designers \citep{hoyle_2014,lewis_2006}.

In RBO some of the constraints become probabilistic, which need specialized methods to handle. The reason is that the probabilistic constraints are expressed by multiple integrals of the joint probability density function (PDF) of random variables, both of which are either practically impossible to obtain or very difficult to evaluate \citep{achintya_2000}. Interested readers are advised to refer to \cite{valdebenito_survey_2010} for a complete review of RBO methods to overcome such difficulties. Within the scope of this paper, we only focus on the first-order reliability methods (FORM), the second-order reliability methods (SORM), the SORA method, and the stochastic response surface (SRS) method. FORM appeared early \citep{cornell_1969} together with the concept of reliability index \citep{hasofer_1974} to solve RBO problems. SORM \citep{fiessler_1979} followed to improve accuracy of the FORM in case of highly nonlinear limit state functions and/or slow decay of the joint PDF. The main idea of FORM and SORM is to approximate the limit state functions using first-order and second-order Taylor series, respectively, at appropriate values (e.g., means) of random variables. This results in a double-loop optimization problem to find the most probable point (MPP). In the context of RBTO, directly solving the double-loop optimization problem has been shown in \cite{maute_reliability-based_2003} for MEMS mechanisms with stochastic loading, boundary conditions as well as material properties; in \cite{sato_reliability-based_2018} for shape uncertainty; in \cite{kang_reliability-based_2018} for geometric imperfections; in \cite{mogami_reliability-based_2006} for frame structures using system reliability under random-variable inputs; in \cite{kang_reliabilitybased_2004} for electromagnetic systems; and in \cite{jung_reliability-based_2004} for geometrically nonlinear structures. This approach is prohibitively expensive and lacks robustness when a large number of random variables present \citep{schueller_critical_2004}. For this reason, single-loop approaches have been developed, in which, e.g., the Karush$-$Kuhn$-$Tucker (KKT) optimality conditions are utilized to avoid the inner loop. Both \cite{nguyen_single-loop_2011} and \cite{silva_component_2010} used variants of the single-loop method in \cite{liang_single-loop_2004} for component and system  reliability-based TO. \cite{kharmanda_reliability-based_2004} is somewhat unique when using their own single-loop method \citep{kharmanda_efficient_2002}. \cite{kogiso_reliability-based_2010} applied the single-loop-single-vector method \citep{chen_reliability_1997} for frame structures under random-variable loads and nonstrucutral mass. Another way to bypass the double-loop problem is the decoupling approaches \citep{valdebenito_survey_2010}, in which reliability analysis results are used to facilitate the optimization loops. Among them, the SORA method is known for its simple implementation compared to the above single-loop methods, and efficiency with FORM \citep{du_sequential_2004,lopez_reliability-based_2012}. This method was employed for RBTO under random-variable inputs in \cite{zhao_reliability-based_2015} and \cite{zhao_comparison_2016}. Simulation techniques coupled with meta modeling or surrogate modeling have enjoyed considerable popularity within RBO community \citep{wang_review_2007}, but has received little attention in TO literature. In \cite{patel_classification_2012}, reliability was assessed using probabilistic neural network classifier for truss structures under random Young's modulus.

To the best of our knowledge, \cite{zhao_reliability-based_2015}, \cite{jalalpour_efficient_2016}, and \cite{keshavarzzadeh_topology_2017} are three papers closest to ours. However, random-field material model was not considered in \cite{zhao_reliability-based_2015} and \cite{keshavarzzadeh_topology_2017}. Furthermore, several concerns can be identified from \cite{keshavarzzadeh_topology_2017}. One of the most important stages in their method is the approximations of failure probability and its sensitivity, which depend on Monte Carlo sampling, and the value of the parameter $\epsilon$ when replacing the Heaviside function with a smooth approximation \citep{keshavarzzadeh_gradient_2016}. Direct Monte Carlo sampling is well-known to have variability \citep{taflanidis_efficient_2008}, meaning two independent runs are very likely to get different values of failure probability and its sensitivity which would obviously affect the optimization results. The parameter $\epsilon$ was chosen by a ``recommendation" backed by observation only. \cite{jalalpour_efficient_2016} assumed that the modeling random field has known marginal distribution, and random variability of Young's modulus is small in order to apply perturbation technique. Both of these assumptions clearly restrict the general applicability of their method. As described in the following sections, our proposed method considers random field uncertainty with the Karhunen$-–$Lo\`{e}ve (KL) expansion used to reduce the dimension of the random field. The KL expansion covers a large class of random field without any restrictions on random variability. In this way, we are able to use the FORM-based inverse reliability method within the SORA framework coupled with the stochastic response surface method to avoid the aforementioned drawbacks of direct Monte Carlo sampling.
\section{Topology optimization under uncertainty}
\label{tou}
\subsection{Deterministic topology optimization}
\label{deter}
This paper adopts a standard notation, which denotes matrices and vectors as bold upper and lower case letters respectively. The below formulation shows a density-based deterministic topology optimization:
\begin{equation}
    \begin{aligned}
        \min_{\boldsymbol{\rho}} && &V( \boldsymbol{\rho} ) = \boldsymbol{\nu}^T \boldsymbol{\rho} \\
        \text{subject to} && &\text{\bf{K}}(\boldsymbol{\rho}) \text{\bf{u}}(\boldsymbol{\rho}) = \text{\bf{f}}, \\
        && &u_i( \boldsymbol{\rho} ) \leq u^0_i, \; i = 1, 2, \dots, m\\
        && &0 < \rho_{\text{min}} \leq \boldsymbol{\rho} \leq 1.
        \label{e:1}
    \end{aligned}
\end{equation}
where $\boldsymbol{\rho}$ is the vector of deterministic finite-element densities; $V( \boldsymbol{\rho} )$ is the total volume of the finite-element mesh; $\boldsymbol{\nu}$ is the vector of finite-element volumes for a unit density; $\text{\bf{K}}(\boldsymbol{\rho})$, $\text{\bf{u}}(\boldsymbol{\rho})$, and $\text{\bf{f}}$ are the stiffness matrix, the displacement vector, and the external load vector, respectively; $u_i( \boldsymbol{\rho} )$ and $u^0_i$, $ i = 1, 2, \dots, m$, are the actual displacement and the maximum allowable displacement at the $i^{th}$ degree of freedom. In the optimization problem~(\ref{e:1}),  $\boldsymbol{\rho}$ is the design variables and $V( \boldsymbol{\rho} )$ is the objective function. The first constraint expresses the equilibrium of the structure while the third constraint is a component-wise inequality, in which each density (design variable)  must be between 1 and a lower limit (e.g., $\rho_{\text{min}} = 0.001$).

To promote manufacturing of the resulting structure, a black-and-white design is preferred. Hence, the Solid Isotropic Material with Penalization (SIMP) method \citep{bendsoe_topology_2003} is used to penalize the intermediate values. According to SIMP, the Young's modulus of each finite elements $E_i$ is obtained as $E_i = \rho_i^pE_i^0$, where $p$ is the penalization factor and $E_i^0$ is the initial value of the Young's modulus corresponding to unit density. The possible values of $p$ were explained in \cite{bendsoe_material_1999}. Applying SIMP transforms the problem~(\ref{e:1}) into a nonlinear optimization problem, which can be solved effectively by many gradient-based algorithms. The Method of Moving Asymptotes (MMA) \citep{svanberg_method_1987,svanberg_class_2002} is chosen in this paper due to its reliability to catch extremum in various settings of TO. Known problems with SIMP are checkerboarding, mesh dependence, and local minima \citep{sigmund_numerical_1998}. Checkerboarding and mesh dependence can be prevented
by mesh independent filtering methods \citep{sigmund_morphology-based_2007}. This paper uses the density filtering \citep{bruns_topology_2001,bourdin_filters_2001} as implemented in \cite{andreassen_efficient_2011}. In the next section we will discuss the changes needed to integrate uncertainty into the optimization problem.
\subsection{Reliability-based topology optimization}
\label{relia}
When uncertainty is introduced into the problem~(\ref{e:1}) in the form of a modeling random field $y(\omega,\mathbf{x})$, the second constraints have to be replaced as follows:
\begin{equation}
    \begin{aligned}
        \min_{\boldsymbol{\rho}} && &V( \boldsymbol{\rho} ) = \boldsymbol{\nu}^T \boldsymbol{\rho} \\
        \text{s.t.} && &\text{\bf{K}}(\boldsymbol{\rho},y) \text{\bf{u}}(\boldsymbol{\rho},y) = \text{\bf{f}}, \\
        && &P_i[u^0_i - u_i( \boldsymbol{\rho},y) < 0] \leq P_i^0, \; i = 1, 2\dots, m\\
        && &0 < \rho_{\text{min}} \leq \boldsymbol{\rho} \leq 1.
        \label{e:relia1}
    \end{aligned}
\end{equation}
where $\mathbf{x} \in D \subset \mathds{R}^d$ is coordinates of a point in a $d$-dimensional physical domain $D$, $\omega \in \Omega$ is an element of the sample space $\Omega$, $P_i[\boldsymbol{\cdot}]$ is the probability, $P_i^0$ is the target probability of the $i^{th}$ constraint. The limit state function is defined as $g(\boldsymbol{\rho},y)=u^0 - u( \boldsymbol{\rho},y)$, so $g(\boldsymbol{\rho},y) < 0$ means constraint violation or failure of the structure in~(\ref{e:1}) while $P_i[g_i(\boldsymbol{\rho},y) < 0]$ shows the failure probability of the $i^{th}$ constraint. The target probability $P_i^0$ is the upper bound of the failure probability $P_i$ and often defined as $P_i^0=\Phi(-\beta_i)$, where $\beta_i$ is the reliability index and $\Phi(\cdot)$ is the standard normal cumulative distribution function. 

Several complications arise when a RBO problem involves a random field. In many cases the considered random field is a infinite-dimensional set, which makes it extremely difficult to propagate uncertainty from input to output of a system if used directly. Instead, the random field usually has to be discretized or dimensionally reduced to a manageable order, which we will show in Sect.~\ref{kle}. The double-loop RBO problem can be very expensive to solve, which is considered in Sect.~\ref{irs}. The probabilistic constraint is often implicit and very hard to evaluate because of complex geometry of its domain, which becomes much easier and cheaper to compute if approximated by the stochastic response surface (SRS) method in Sect.~\ref{srs}. This section ends with a thorough description of our proposed algorithm for a RBTO problem based on the above methods.
\subsubsection{Karhunen$-–$Lo\`{e}ve expansion}
\label{kle}
A random field can be discretized by serveral methods including the Expansion Optimal Linear Estimator \citep{li_optimal_1993} and polynomial chaos expansion \citep{xiu2010numerical,ghanem2003stochastic}. Compared to others, the Karhunen$-–$Lo\`{e}ve (KL) expansion \citep{loeve2017probability} is ``the most effcient in terms of the number of random variables required for a given accuracy" \citep{sudret_2000}. The KL expansion of a random field $y(\omega,\mathbf{x})$ is given as
\begin{equation}
    y(\omega,\mathbf{x}) = E[\mathbf{x}] + \sum_{i=1}^\infty\sqrt{\lambda_i}\xi_i(\omega)e_i(\mathbf{x})
    \label{e:2}
\end{equation}
where $E[\mathbf{x}]$ is the mean of the random field. The orthogonal eigenfunctions $e_i(\mathbf{x})$ and the corresponding eigenvalues $\lambda_i$ are solutions of the  following eigenvalue problem:
\begin{equation}
    \int_{D}K(\mathbf{x}_1,\mathbf{x}_2)e_i(\mathbf{x})d\mathbf{x}=\lambda_ie_i(\mathbf{x}) \quad \mathbf{x},\mathbf{x}_1,\mathbf{x}_2 \in D
    \label{e:3}
\end{equation}
where $K(\mathbf{x}_1,\mathbf{x}_2)$ is the covariance function of the random field
\begin{equation}
    K(\mathbf{x}_1,\mathbf{x}_2) = E\left[y(\mathbf{x}_1)y(\mathbf{x}_2)\right] \quad \mathbf{x}_1,\mathbf{x}_2 \in D
\end{equation}
The random variables $\xi_i(\omega)$ are uncorrelated and satisfy:
\begin{equation}
    \begin{aligned}
        & E[\xi_i] = 0, E[\xi_i\xi_j] = \delta_{ij} \\
        & \xi_i(\omega) = \frac{1}{\sqrt{\lambda_i}}\int_D\left(y(\omega,\mathbf{x})-E[\mathbf{x}]\right)e_i(\mathbf{x})d\mathbf{x}
    \end{aligned}
\end{equation}
where $\delta_{ij}$ is the  Kronecker delta. The infinite series in~(\ref{e:2}) has to be truncated to use in practice. Based on the fact that the influence of higher order terms decays rapidly, only a few of the terms are needed to capture behavior of the random field with appropriate precision.

The KL expansion requires the solution of the eigenvalue problem~(\ref{e:3}), which is pretty straightforward in the case of a random process (1-dimensional random field) \citep{ray_numerical_2013}. For the purpose of demonstration and without loss of generality, this paper in Sect.~\ref{results} exploits the separability of the covariance function of a 2-dimensional random field:
\begin{equation}
    \begin{aligned}
    K(\mathbf{s},\mathbf{t}) &= exp\left(\frac{-|s_1-t_1|}{l_1}\times\frac{-|s_2-t_2|}{l_2}\right) \\
    &=exp\left(-\frac{|s_1-t_1|}{l_1}\right)exp\left(-\frac{|s_2-t_2|}{l_2}\right)\\
    &\mathbf{s},\mathbf{t} \in D \subset \mathds{R}^2
    \label{e:covariance}
    \end{aligned}
\end{equation}
where $l_1$ and $l_2$ are the correlation lengths in the two coordinate directions. This class of covariance function makes the eigenvalues and eigenfunctions separable in the sense that both are the product of their univariate counterparts \citep{wang_2008}.

\subsubsection{Inverse reliability and SORA}
\label{irs}
\begin{figure}[ht]
    \centering
    \includegraphics[scale=0.85]{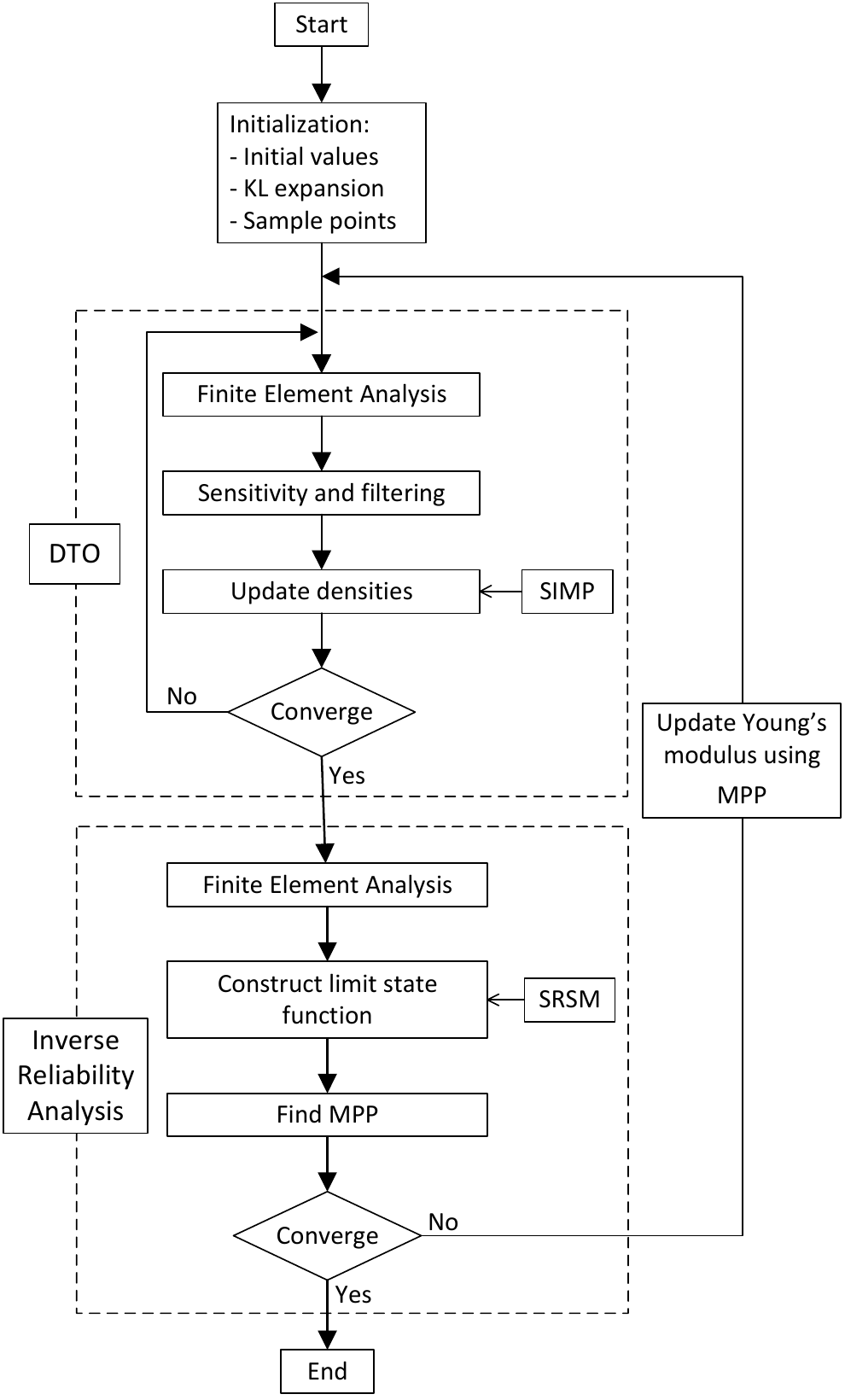}
    \caption{SORA-based RBTO flowchart \citep{zhao_reliability-based_2015}.}
    \label{fig:sora}
\end{figure}

The form of the probabilistic constraints in~(\ref{e:relia1}) is called the reliability index approach (RIA). However, according to \cite{tu_new_1999}, the performance measure approach (PMA) provides better numerical stability and higher rate of convergence. Using the PMA, (\ref{e:relia1}) is rewritten as
\begin{equation}
    \begin{aligned}
        \min_{\boldsymbol{\rho}} && &V( \boldsymbol{\rho} ) = \boldsymbol{\nu}^T \boldsymbol{\rho} \\
        \text{s.t.} && &\text{\bf{K}}(\boldsymbol{\rho},y) \text{\bf{u}}(\boldsymbol{\rho},y) = \text{\bf{f}}, \\
        && &u^0_i - u_i( \boldsymbol{\rho},y) \geq 0, \; i = 1, 2\dots, m\\
        && &0 < \rho_{\text{min}} \leq \boldsymbol{\rho} \leq 1.
        \label{e:relia2}
    \end{aligned}
\end{equation}
Solving the above problem requires a truncated KL expansion $y(\omega,\mathbf{x}) \approx y(\xi_i(\omega),\mathbf{x})$, FORM, and inverse reliability analysis. In order to apply FORM, the random vector $\mathbf{\Xi}=\{\xi_i\}$ is transformed into a vector of standard normal random variables $\mathbf{\Psi}=\{\psi_i\}$ using the Rosenblatt or the Nataf transformation $\mathbf{\Psi} = T(\mathbf{\Xi})$ or $\mathbf{\Xi} = T^{-1}(\mathbf{\Psi})$. Then, the most probable point (MPP) $\boldsymbol{\xi}_i^*$ in physical space or $\boldsymbol{\psi}_i^*$ in transformed space is obtained by inverse reliability analysis:
\begin{equation}
    \begin{aligned}
        \min_{\boldsymbol{\psi}} && &g_i( \boldsymbol{\psi} ) \\
        \text{s.t.} && &\parallel\boldsymbol{\psi}\parallel=\beta_i.
        \label{e:relia3}
    \end{aligned}
\end{equation}
where $g_i( \boldsymbol{\psi} )$ is the $i^{th}$ limit state function. In this paper, the Matlab CODES toolbox \citep{codes_toolbox} is chosen to solve~(\ref{e:relia3}). Also, the SORA framework is adopted to decouple the double-loop structure of~(\ref{e:relia2}). In SORA, instead of nesting the optimization problem~(\ref{e:relia3}) within~(\ref{e:relia2}), it serializes~(\ref{e:relia2}) into a chain of loops of deterministic TO and inverse reliability analysis (Fig.~\ref{fig:sora}). Each $k^{th}$ loop starts with deterministic TO followed by inverse reliability analysis:
\begin{equation}
    \begin{aligned}
        \min_{\boldsymbol{\rho}^k} && &V( \boldsymbol{\rho}^k ) = \boldsymbol{\nu}^T \boldsymbol{\rho}^k \\
        \text{s.t.} && &\text{\bf{K}}(\boldsymbol{\rho},y) \text{\bf{u}}(\boldsymbol{\rho},y) = \text{\bf{f}}, \\
        && &u^0_i - u_i\left( \boldsymbol{\rho}^k,y(\boldsymbol{\xi}_i^{*(k-1)},\mathbf{x})\right) \geq 0, \; i = 1, 2\dots, m\\
        && &0 < \rho_{\text{min}} \leq \boldsymbol{\rho}^k \leq 1.
        \label{e:relia4}
    \end{aligned}
\end{equation}
where $\boldsymbol{\xi}_i^{*(k-1)}$ denotes the MPP in physical space of $i^{th}$ limit state function in the $(k-1)^{th}$ loop. Solving~(\ref{e:relia4}) gives $\boldsymbol{\rho}^{*(k)}$, which is substituted into~(\ref{e:relia3}) to find the next MPP $\boldsymbol{\xi}_i^{*(k)}$ in the form of $\boldsymbol{\psi}_i^{*(k)}$:
\begin{equation}
    \begin{aligned}
        \min_{\boldsymbol{\psi}} && &g_i(\boldsymbol{\rho}^{*(k)},\boldsymbol{\psi} ) \\
        \text{s.t.} && &\parallel\boldsymbol{\psi}\parallel=\beta_i.
        \label{e:relia5}
    \end{aligned}
\end{equation}
The computational cost is saved by a significant reduction in the number of reliability analyses needed until both the deterministic TO and reliability analysis converge. One obvious problem of the above scheme is that the final $\boldsymbol{\rho}^*$ may not be the same as in the double-loop setting, which is compensated for by shifting the random design variables in each loop~\citep{yin_enhanced_2006}. Since there is no random design variable (only random parameters) in the setting of this paper, no such modification is needed.
\subsubsection{Stochastic response surface method}
\label{srs}
Considering random inputs, the propagated uncertainty in many cases makes it difficult to get an explicit expression of interested outputs, which blocks further progress, such as taking derivatives and optimization. The SRS method overcomes this problem by approximating the outputs by a polynomial chaos expansion \citep{xiu2010numerical}. The below formulations follows \cite{huang_extended_2007}. The multidimensional Hermite polynomials of degree $p$ are used in the SRS method and defined as:
\begin{equation}
    \begin{aligned}
    &H_p(\alpha_{i_1},\alpha_{i_2},\ldots,\alpha_{i_p})\\
    &=(-1)^pe^{\frac{1}{2}\boldsymbol{\alpha}^T\boldsymbol{\alpha}}\frac{\partial^p}{\partial\alpha_{i_1},\partial\alpha_{i_2},\ldots,\partial\alpha_{i_p}}e^{-\frac{1}{2}\boldsymbol{\alpha}^T\boldsymbol{\alpha}}
    \end{aligned}
    \label{e:hermite}
\end{equation}
where $\boldsymbol{\alpha}=\{\alpha_{i_k}\}_{k=1}^p$ is a vector of standard normal random variables. The interested output $z$ is estimated as follows:
\begin{equation}
    \begin{aligned}
    z = a_0 &+ \sum_{i_1=1}^n a_{i_1}H_1(\alpha_{i_1}) + \sum_{i_1=1}^n\sum_{i_2=1}^{i_1} a_{i_1i_2}H_2(\alpha_{i_1},\alpha_{i_2})\\
    &+\sum_{i_1=1}^n\sum_{i_2=1}^{i_1}\sum_{i_3=1}^{i_2} a_{i_1i_2i_3}H_3(\alpha_{i_1},\alpha_{i_2},\alpha_{i_3}) + \ldots
    \end{aligned}
    \label{e:srs}
\end{equation}
where $n$ is the number of standard normal random variables used in the expansion, and $a_0,a_{i_1},a_{i_1i_2},a_{i_1i_2i_3},\ldots$ are unknown coefficients. If $n=2$ and $p=3$, then the expansion~(\ref{e:srs}) will become:
\begin{equation}
    \begin{aligned}
    z(\alpha_{i_1},\alpha_{i_2})&=&a_0&+a_1\alpha_{i_1}+a_2\alpha_{i_2}+a_3(\alpha_{i_1}^2-1)+a_4(\alpha_{i_2}^2-1)\\
    && &+a_5\alpha_{i_1}\alpha_{i_2}+a_6(\alpha_{i_1}^3-3\alpha_{i_1})+a_7(\alpha_{i_2}^3-3\alpha_{i_2})\\
    && &+a_8(\alpha_{i_1}\alpha_{i_2}^2-\alpha_{i_1})+a_9(\alpha_{i_1}^2\alpha_{i_2}-\alpha_{i_2})\\
    &=&a_0&+\sum_{k=1}^9 a_k\eta_{i_k}
    \end{aligned}
    \label{e:n2p3}
\end{equation}
where $1,\eta_{i_1},\eta_{i_2},\ldots,\eta_{i_9}$ are Hermite polynomials. The ten unknown coefficients $a_0,a_1,\ldots,a_9$ are found by solving a system of linear equations using at least ten different realizations of $(\alpha_{i_1},\alpha_{i_2})$. A stochastic response surface constructed with 17 collocation points was shown to be a very good approximation in \cite{huang_extended_2007}, which will be used in this paper.

\subsubsection{Solution algorithm}
\label{algorithm}
A step-by-step explanation of our solution algorithm (Fig.~\ref{fig:sora}) is listed below:
\begin{enumerate}
    \item Initialize the problem: finite element mesh; initial values of design variables, SIMP and optimization parameters; Karhunen$-–$Lo\`{e}ve (KL) expansion of random field; sample points for the stochastic response surface method (SRSM); etc.
    \item Deterministic topology optimization (DTO): the SIMP method and the MMA are used to solve~(\ref{e:relia4}), which is deterministic because there is no random variable involved. To guarantee that the Young's modulus is physically meaningful (e.g., only positive values), it is modeled as in \cite{lazarov_topology_2012}:
    \begin{equation}
    E(\mathbf{x}) = F^{-1} \circ \Phi\left[y(\omega,\mathbf{x})\right]
    \label{e:young}
    \end{equation}
    where $\Phi[\cdot]$ is the standard normal cumulative distribution function (CDF) and $F^{-1}$ is the inverse of a prescribed CDF. The uniform distribution is chosen in this paper resulting in
    \begin{equation}
    E(\mathbf{x}) = a + (b-a)\Phi\left[y(\omega,\mathbf{x})\right]
    \label{e:uniform}
    \end{equation}
    where $a$ and $b$ are the two bounds of the distribution. Two other physically admissible distributions are the log-normal and the beta distribution, which can replace the uniform distribution in~(\ref{e:young}) with trivial effort.
    \item Inverse reliability analysis (IRA): the optimum values of design variables (densities) found in the previous stage and the sample points are used to construct response surfaces of the probabilistic constraints, which in turn are utilized in~(\ref{e:relia5}) to find the most probable point (MPP). Based on convergence conditions, the algorithm may stop or a new loop is requested with updated Young's modulus using the MPP.
\end{enumerate}

\section{Results}
\label{results}
In this section two common benchmark problems (the MBB and the L-shaped beam) are used to test our proposed algorithm with three values of target reliability and four different parameter tuples of the uniform distribution in~(\ref{e:uniform}). The optimization results are then verified by Monte Carlo simulations. For simplicity, all quantities are
given dimensionless.

The two numerical examples below share some common implementation settings. Square, linear plane stress element is employed in both examples, which has unit side length and thickness, and is made of isotropic linear elastic material with Poisson's ratio $\nu=0.3$. The material Young's modulus is assumed to be a  centered mean-square Gaussian random field with known covariance function as in~(\ref{e:covariance}). Under such assumption, the KL expansion results in a series of independent standard normal random variables \citep{noauthor_brief_nodate}. The correlation length are chosen as $l_1=l_2=0.6$, and the truncated KL expansion contains the first two eigenvalues and eigenfunctions. The CODES toolbox provides a wide range of optimization algorithms, among which the Hybrid Mean Value (HMV) method is chosen due to its efficiency \citep{youn_hybrid_2003}. The three values of target reliability are $\{2.0,2.5,3.0\}$ while the four parameter tuples are $(a,b) = \{(1,1.1),(1,1.3),(1,1.5),(1,1.7)\}$. The MMA is the optimizer running on design variables in DTO until convergence criterion, e.g., maximum difference of design variables of two consecutive iterations $(d_{\text{max}}=0.001)$, is satisfied. The IRA stage is considered converged if maximum difference of the MPPs of the previous and current loop is smaller than $0.001$, and the number of iterations is smaller then $20$. In SIMP, the minimum length scale $r_{\text{min}}=1.5$ and penalization factor $p=3$ are employed. For verification, $50000$ Monte Carlo simulations (MCS) are run to calculate the failure probabilities and the statistical moments, which are compared with the expected failure probabilities computed from the three reliability levels, and the SRSM-based results. 

\subsection{The MBB beam}
\label{mbb}
\begin{figure}
    \centering
    \includegraphics[width=0.45\textwidth]{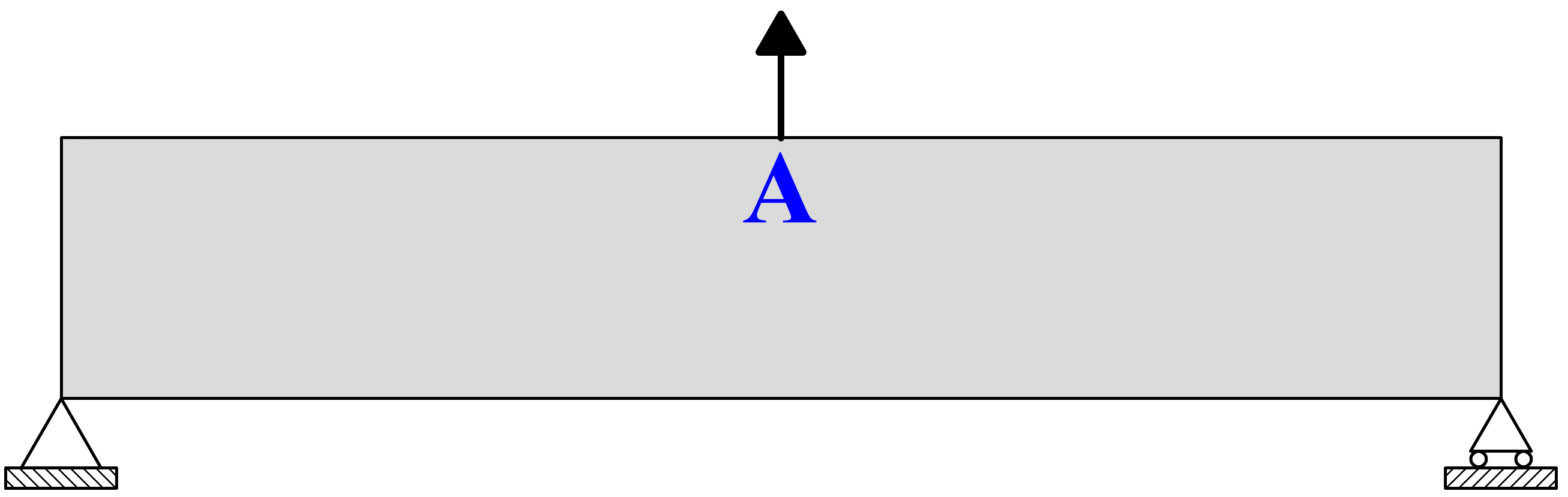}
    \caption{The MBB beam}
    \label{fig:mbb}
\end{figure}

Consider the simply supported, two-dimensional domain in Fig.~\ref{fig:mbb}. It is subject to an unit vertical load at the midpoint of its top edge and meshed into 120 $\times$ 20 elements. For the purpose of this paper, its volume is minimized under a constraint on the displacement of the load application point. The maximum allowable vertical displacement of the load application point (point A in Fig.~\ref{fig:mbb}) is given as $u^0=170$. This value is picked up for demonstration purpose: volume fraction of RBTO results, i.e., ratio of optimized volume to that of initial domain, will be in the range of $(0.4,0.6)$ facilitating visual inspection; any reasonable values should work. Only one half of the domain (60 $\times$ 20 elements) needs to be optimized due to its symmetry.

\begin{sidewaystable*}[ph!]
    \centering
    \caption{The MBB beam: RBTO results}
    \begin{tabular}{lcccc}
        \toprule
        $(a,b)$ & $(1,1.1)$ & $(1,1.3)$ & $(1,1.5)$ & $(1,1.7)$ \\
        \toprule
        $\beta=2$
        & \includegraphics[scale=0.43]{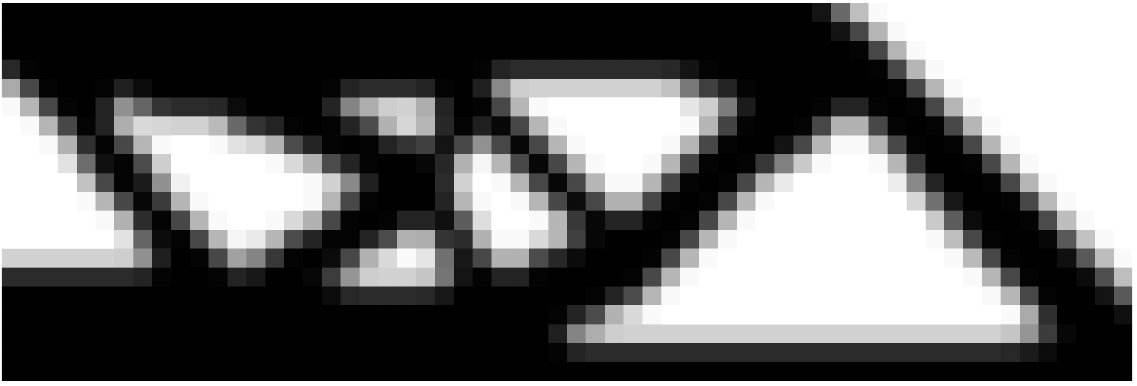}
        & \includegraphics[scale=0.43]{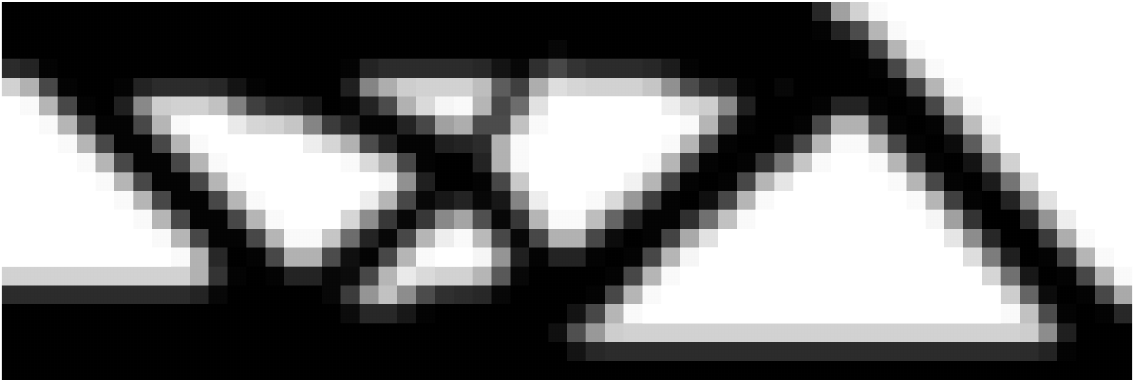}
        & \includegraphics[scale=0.43]{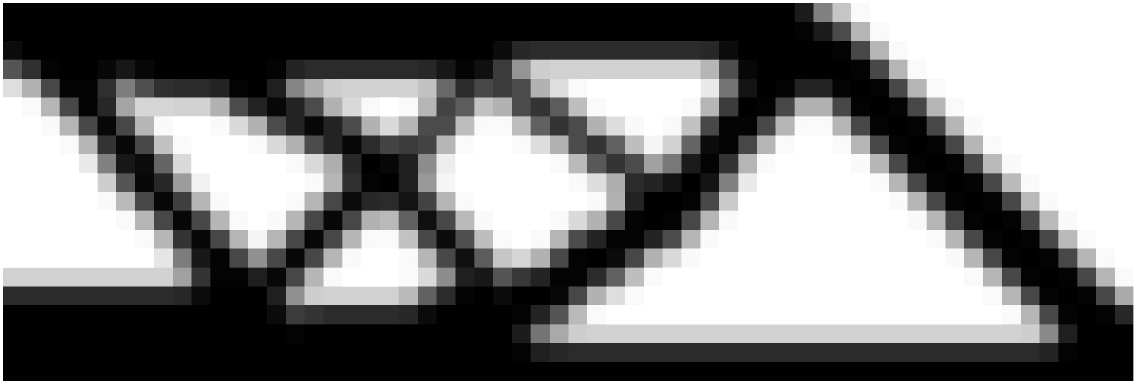}
        & \includegraphics[scale=0.43]{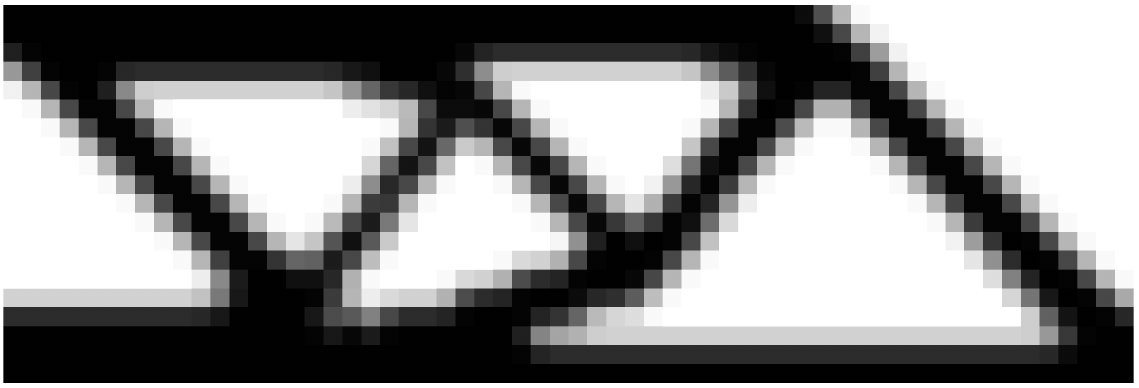} \\
        \midrule
        $\beta=2.5$ 
        & \includegraphics[scale=0.43]{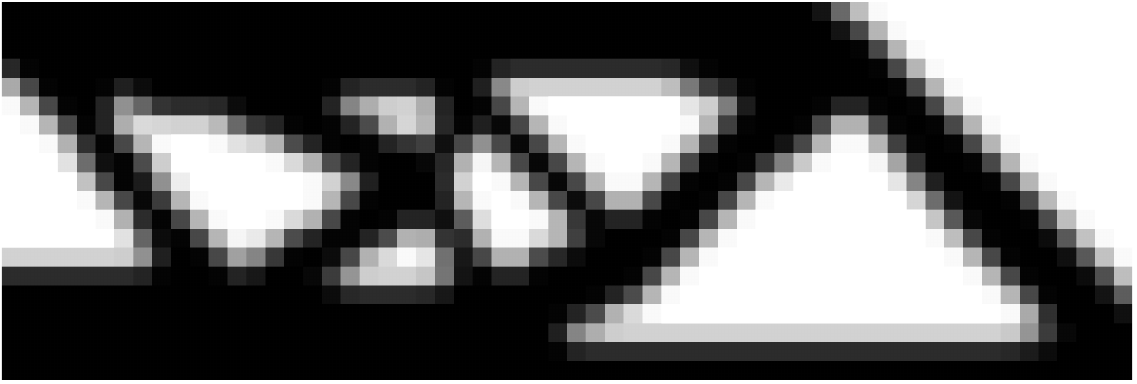}
        & \includegraphics[scale=0.43]{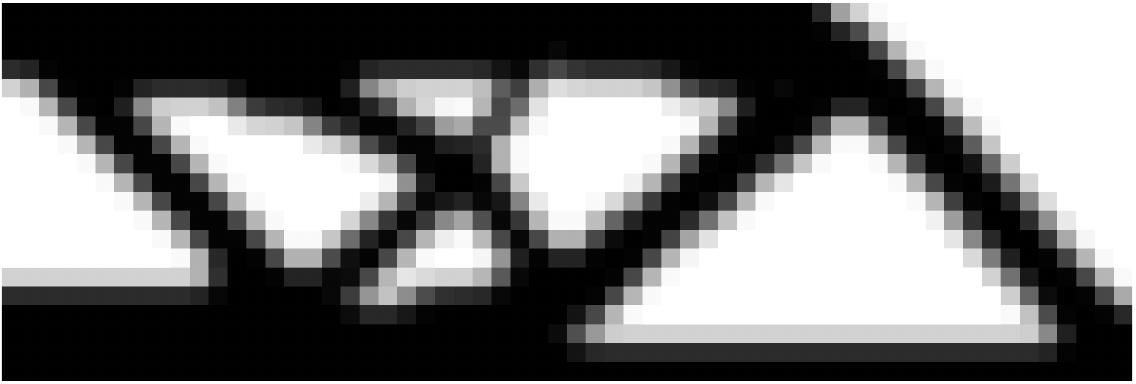} 
        & \includegraphics[scale=0.43]{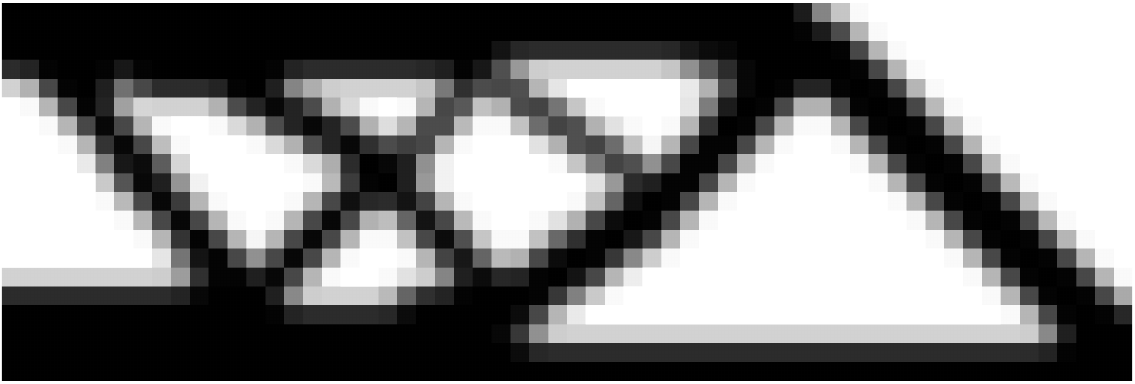}
        & \includegraphics[scale=0.43]{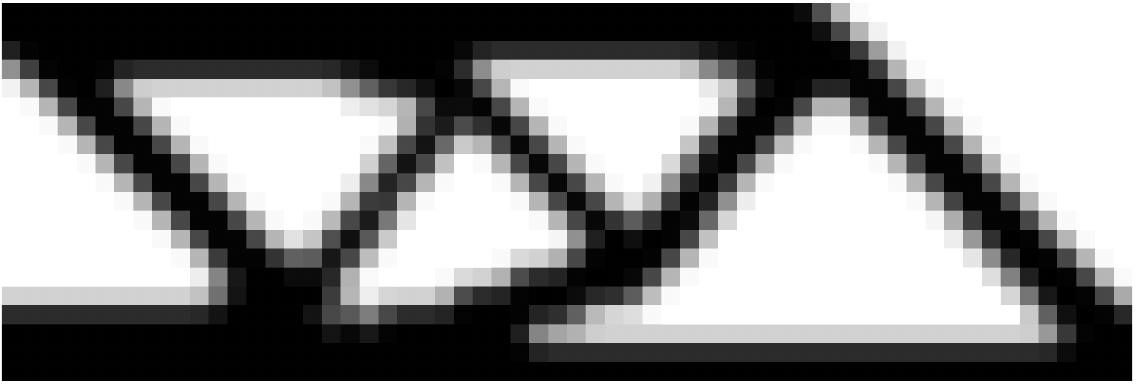} \\
        \midrule
        $\beta=3$ 
        & \includegraphics[scale=0.43]{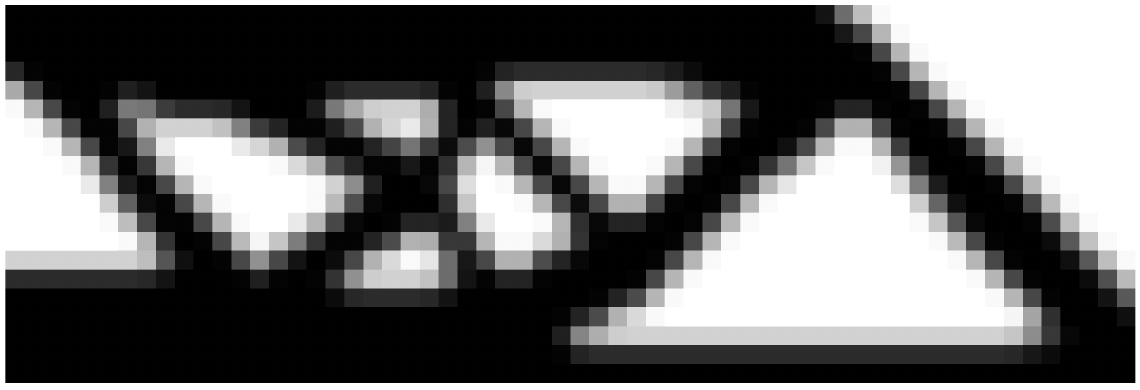}
        & \includegraphics[scale=0.43]{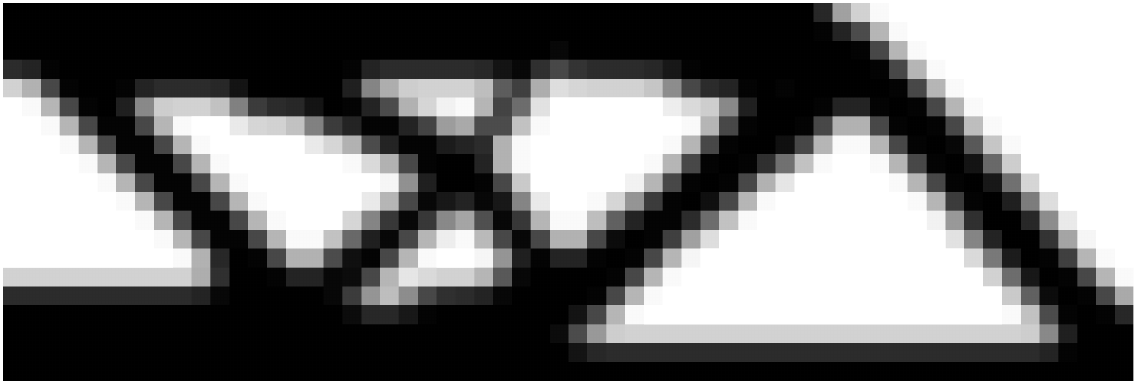} 
        & \includegraphics[scale=0.43]{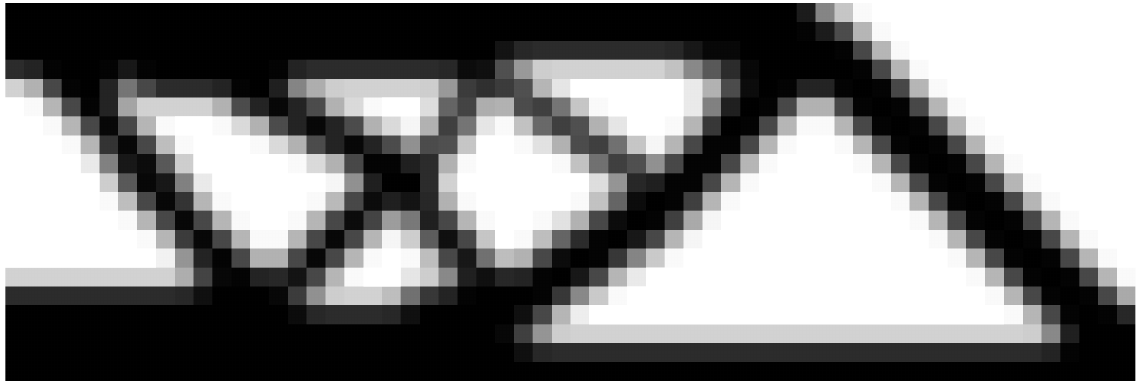}
        & \includegraphics[scale=0.43]{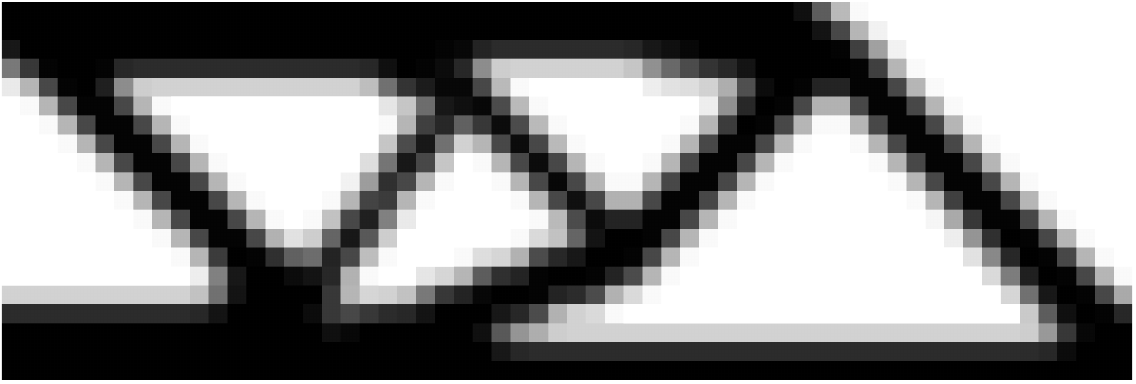} \\
        \bottomrule
        \label{tbl:RBTO_MMB}
    \end{tabular}
    \vspace{10px}
    \caption{The MBB beam: DTO results}
    \begin{tabular}{lcccc}
        \toprule
        $E$  & $1.05$ & $1.15$ & $1.25$ & $1.35$ \\
        \midrule
        DTO 
        & \includegraphics[scale=0.43]{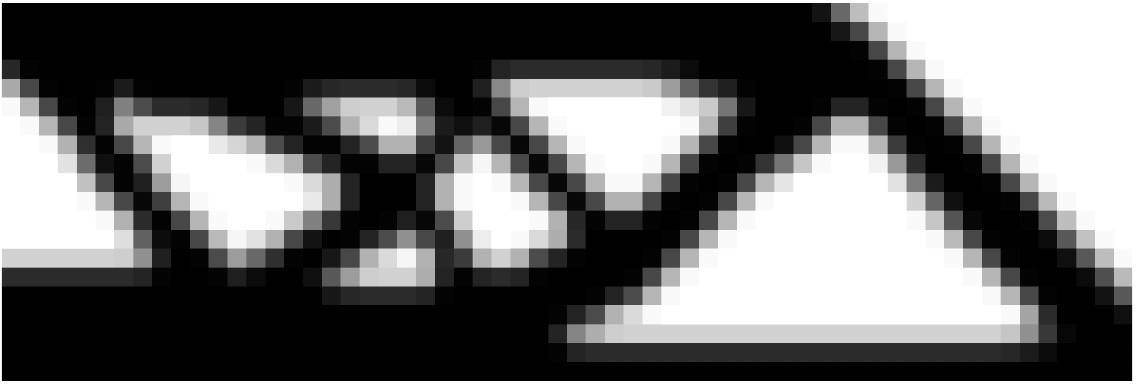}
        & \includegraphics[scale=0.43]{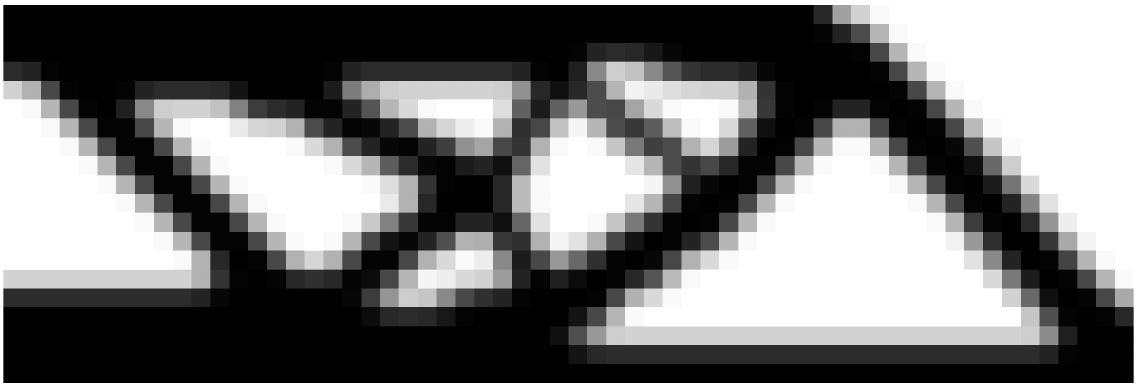}
        & \includegraphics[scale=0.43]{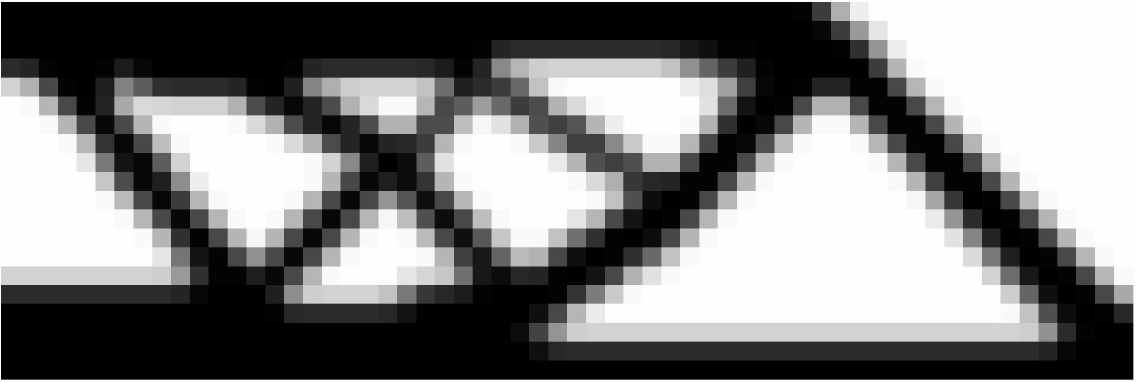}
        & \includegraphics[scale=0.43]{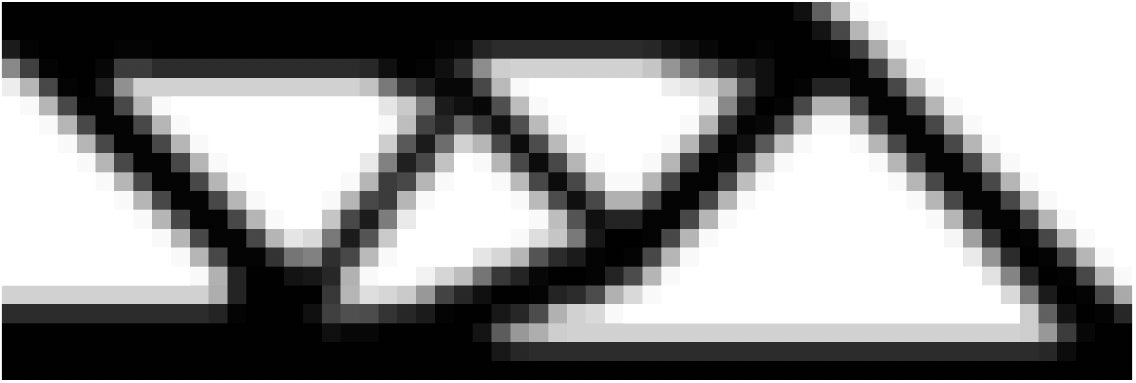} \\
        \midrule
        \multicolumn{1}{m{4em}}{Volume fraction} & 0.5961 & 0.5464 & 0.5051 & 0.4693 \\
        \bottomrule
        \label{tbl:DTO_MMB}
    \end{tabular}
\end{sidewaystable*}

\begin{table*}
    \centering
    \caption{The MBB beam: MCS, SRSM, and volume fraction}
    \begin{tabular}{lllllllll}
        \toprule
        & & & & \multicolumn{3}{c}{MCS} & \multicolumn{2}{c}{SRSM} \\
        \cmidrule(lrr){5-7}
        \cmidrule(lr){8-9}
        $\beta$ & Expected $P_f$ & $(a,b)$ & Volume fraction & $\mu$ & $\sigma$ & $P_f$ & $\mu$ & $\sigma$ \\
        \cmidrule{1-9}
        \multirow{4}{*}[-1em]{2} & \multirow{4}{*}[-1em]{0.02275} & $(1,1.1)$ & 0.5972 & 169.751 & 0.1238 & 0.02192 & 169.751 & 0.1238 \\
        \cmidrule{3-9}
        & & $(1,1.3)$ & 0.5452 & 169.318 & 0.3386 & 0.02242 & 169.318 & 0.3386 \\
        \cmidrule{3-9}
        & & $(1,1.5)$ & 0.5070 & 168.952 & 0.5196 & 0.02242 & 168.952 & 0.5196 \\
        \cmidrule{3-9}
        & & $(1,1.7)$ & 0.4732 & 168.654 & 0.6665 & 0.02276 & 168.654 & 0.6665 \\
        \cmidrule{1-9}
        \multirow{4}{*}[-1em]{2.5} & \multirow{4}{*}[-1em]{0.00620} & $(1,1.1)$ & 0.5974 & 169.689 & 0.1236 & 0.00554 & 169.689 & 0.1236 \\
        \cmidrule{3-9}
        & & $(1,1.3)$ & 0.5457 & 169.147 & 0.3381 & 0.00564 & 169.147 & 0.3381 \\
        \cmidrule{3-9}
        & & $(1,1.5)$ & 0.5074 & 168.682 & 0.5209 & 0.00578 & 168.682 & 0.5209 \\
        \cmidrule{3-9}
        & & $(1,1.7)$ & 0.4740 & 168.315 & 0.6643 & 0.00584 & 168.315 & 0.6643 \\
        \cmidrule{1-9}
        \multirow{4}{*}[-1em]{3} & \multirow{4}{*}[-1em]{0.001349} & $(1,1.1)$ & 0.5977 & 169.628 & 0.1232 & 0.001320 & 169.628 & 0.1232 \\
        \cmidrule{3-9}
        & & $(1,1.3)$ & 0.5463 & 168.976 & 0.3376 & 0.001340 & 168.976 & 0.3376 \\
        \cmidrule{3-9}
        & & $(1,1.5)$ & 0.5081 & 168.411 & 0.5218 & 0.001380 & 168.411 & 0.5218 \\
        \cmidrule{3-9}
        & & $(1,1.7)$ & 0.4742 & 167.970 & 0.6640 & 0.001360 & 167.970 & 0.6640 \\
        \bottomrule
        \label{tbl:MBB_Data}
    \end{tabular}
\end{table*}

For comparison with RBTO results, the deterministic topology optimization (DTO) is performed with the Young's modulus equal to $\{1.05,1.15,1.25,1.35\}$, which are the means of the corresponding uniform distributions in Table~\ref{tbl:RBTO_MMB}. Table~\ref{tbl:DTO_MMB} shows DTO results of the MBB beam found by solving~(\ref{e:1}), and their volume fractions. Next, the proposed algorithm is run with different values of target reliability and material parameters resulting in 12 optimized designs as shown in Table~\ref{tbl:RBTO_MMB}. To verify our proposed algorithm, the results in Table~\ref{tbl:RBTO_MMB} are used in the MCS and SRSM. The mean $\mu$ and the standard deviation $\sigma$ of vertical displacement of point A, and the failure probabilities $P_f$ of the constraint in~(\ref{e:relia1}) are aggregated in Table~\ref{tbl:MBB_Data} from those simulations. Finally, the cumulative distribution function (CDF) of vertical displacement of point A is plotted. Because the CDF plots constructed from the MCS and SRSM are very close to each other for every combination of reliability levels and material parameters and no additional conclusion can be drawn by including all of them, only one CDF plot in the case of $\beta=2.5$ and $(a,b)=(1,1.5)$ is shown in Fig.~\ref{fig:MBB_CDF_all}. Also, the last 10 points of the simulations, which belong to the tail of the CDF where failure probabilities are determined, are displayed in Fig.~\ref{fig:MBB_CDF_10} to aid visual inspection. A number of conclusions and observations from this example will be discussed in Sect.~\ref{disu}.

\subsection{The L-shaped beam}
\label{L-shaped}
\begin{figure}[H]
    \centering
    \includegraphics[width=0.25\textwidth]{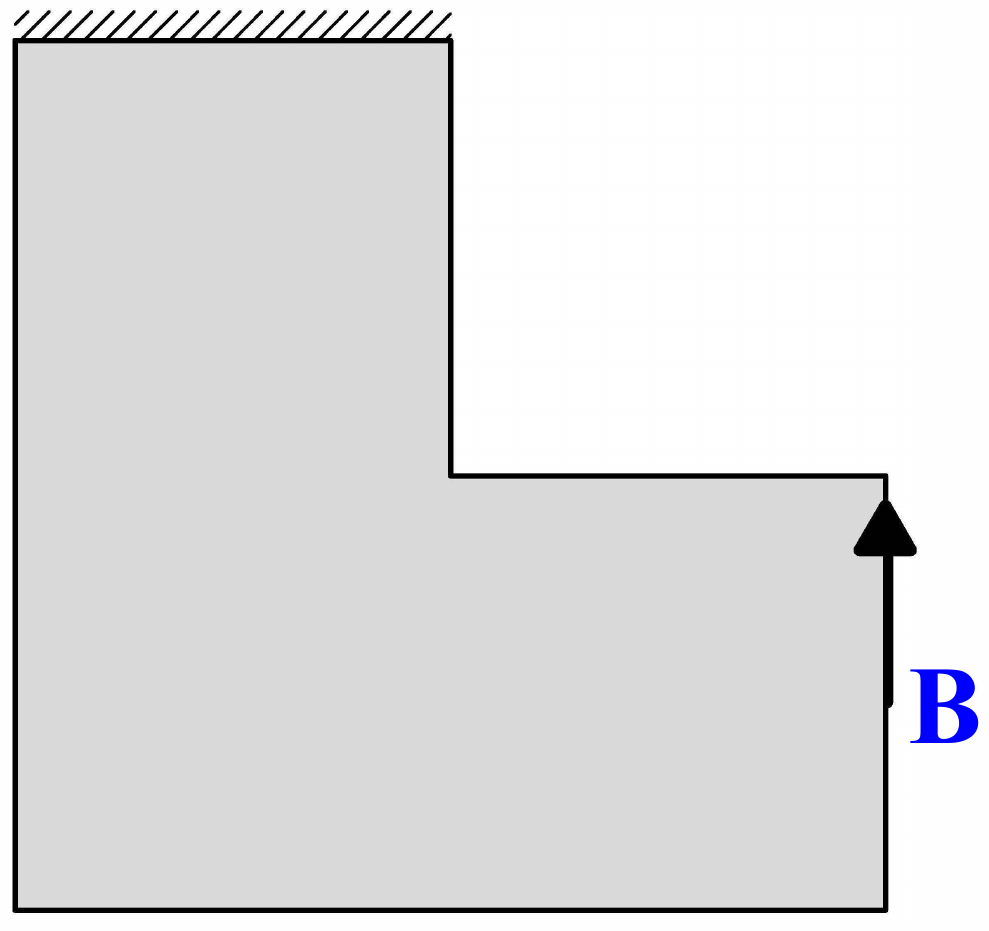}
    \caption{The L-shaped beam}
    \label{fig:L}
\end{figure}
\begin{table*}
    \centering
    \caption{The L-shaped beam: RBTO results}
    \begin{tabular}{lcccc}
        \toprule
        $(a,b)$ & $(1,1.3)$ & $(1,1.5)$ & $(1,1.7)$ \\
        \toprule
        $\beta=2$
        & \includegraphics[scale=0.43]{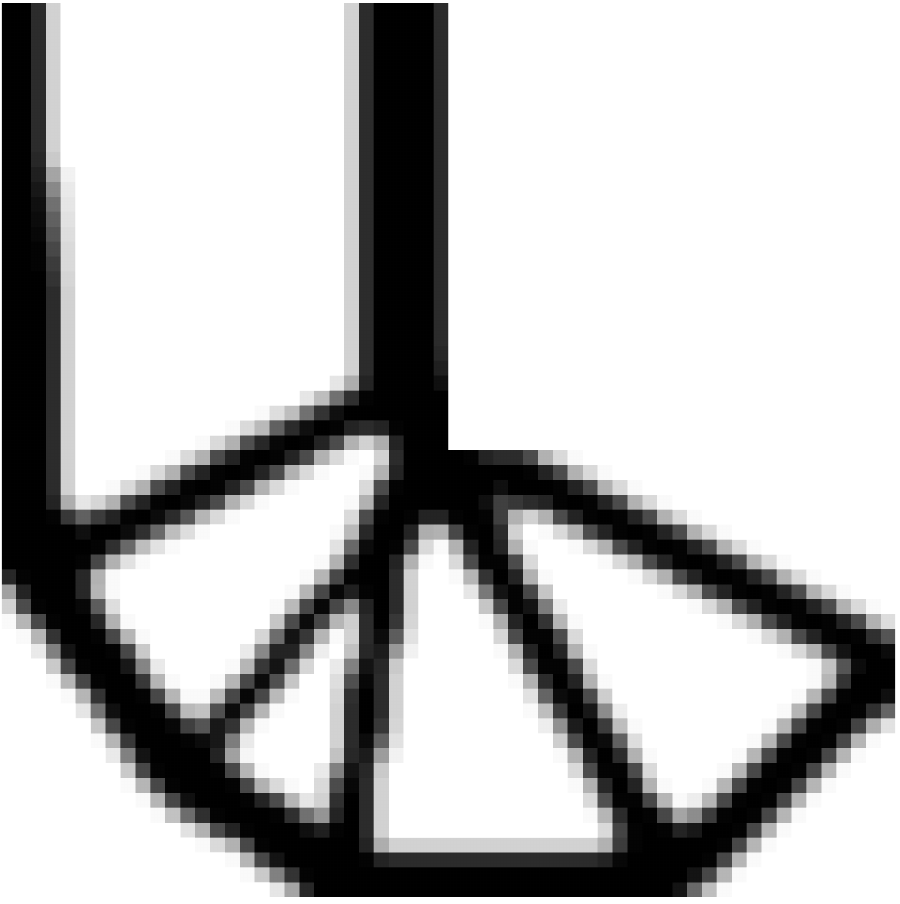}
        & \includegraphics[scale=0.43]{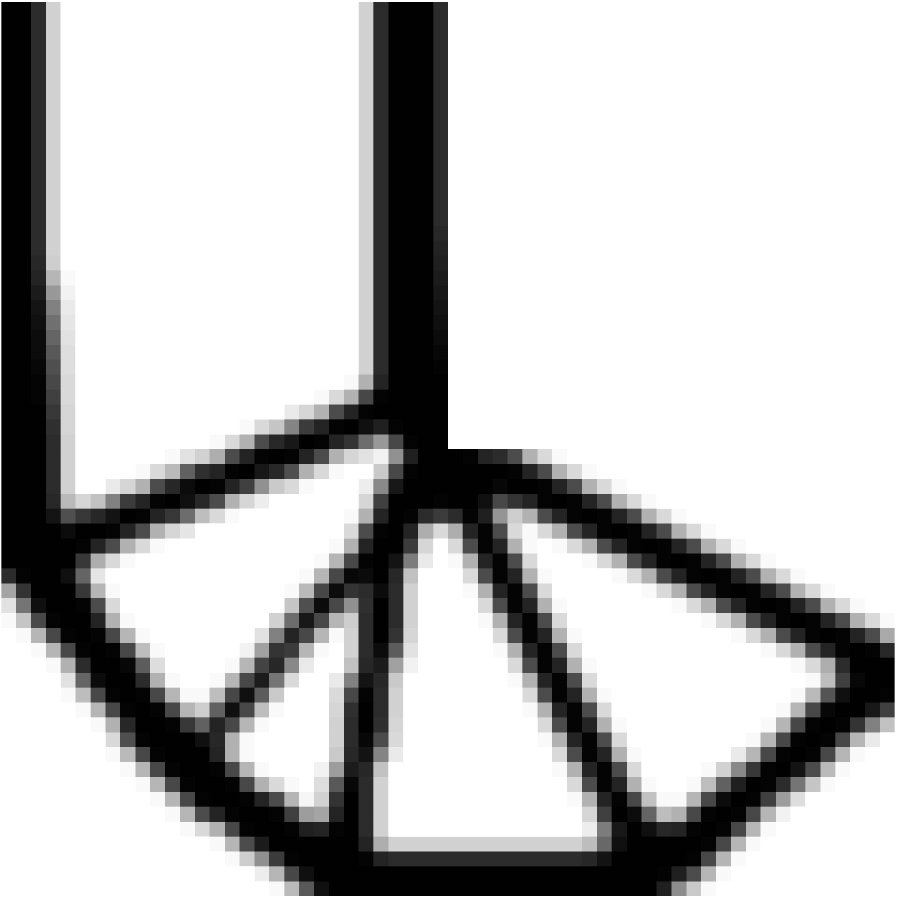}
        & \includegraphics[scale=0.43]{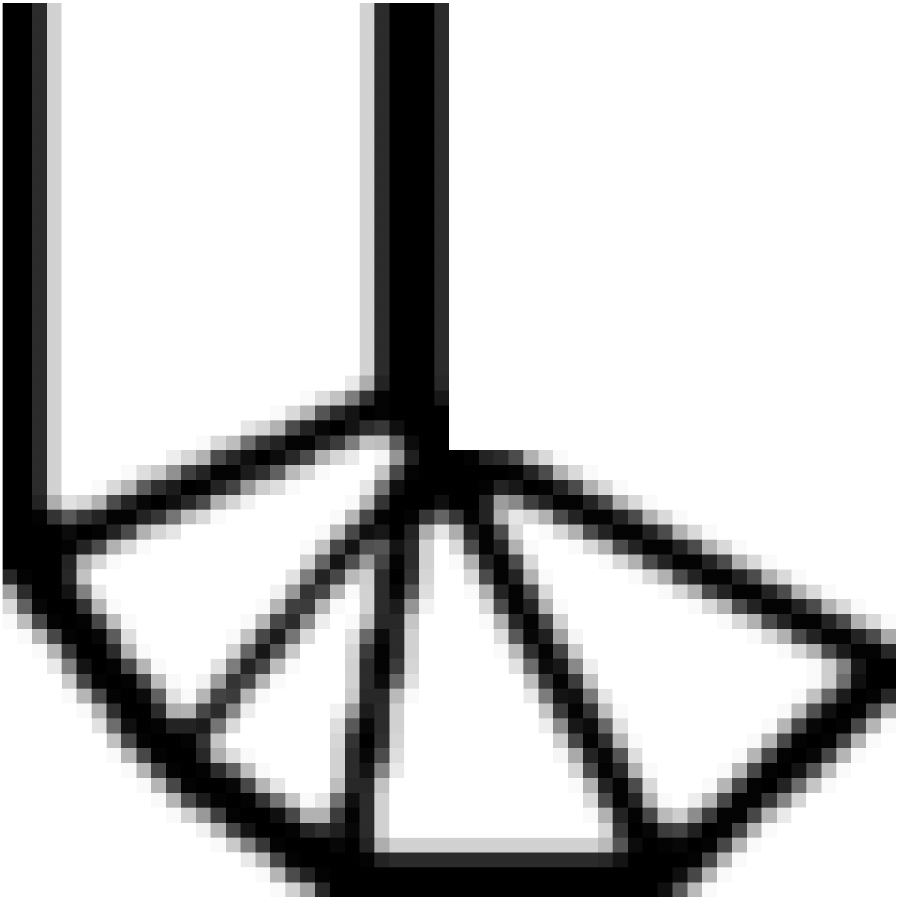} \\
        \midrule
        $\beta=3$ 
        & \includegraphics[scale=0.43]{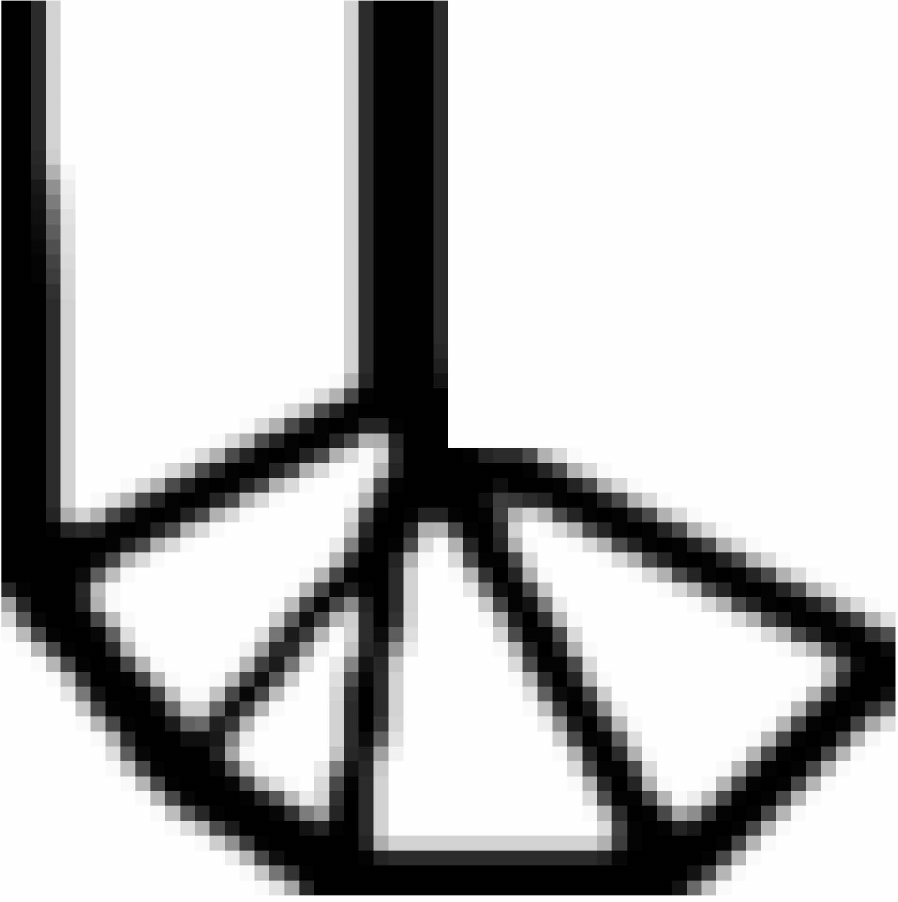} 
        & \includegraphics[scale=0.43]{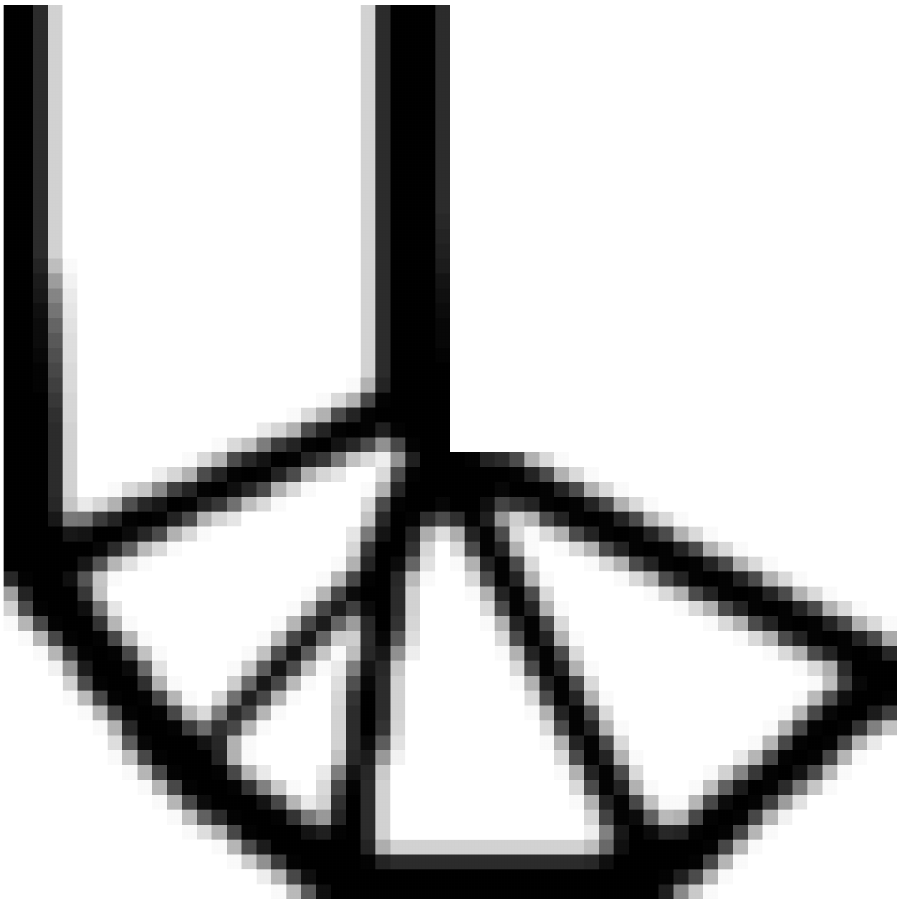}
        & \includegraphics[scale=0.43]{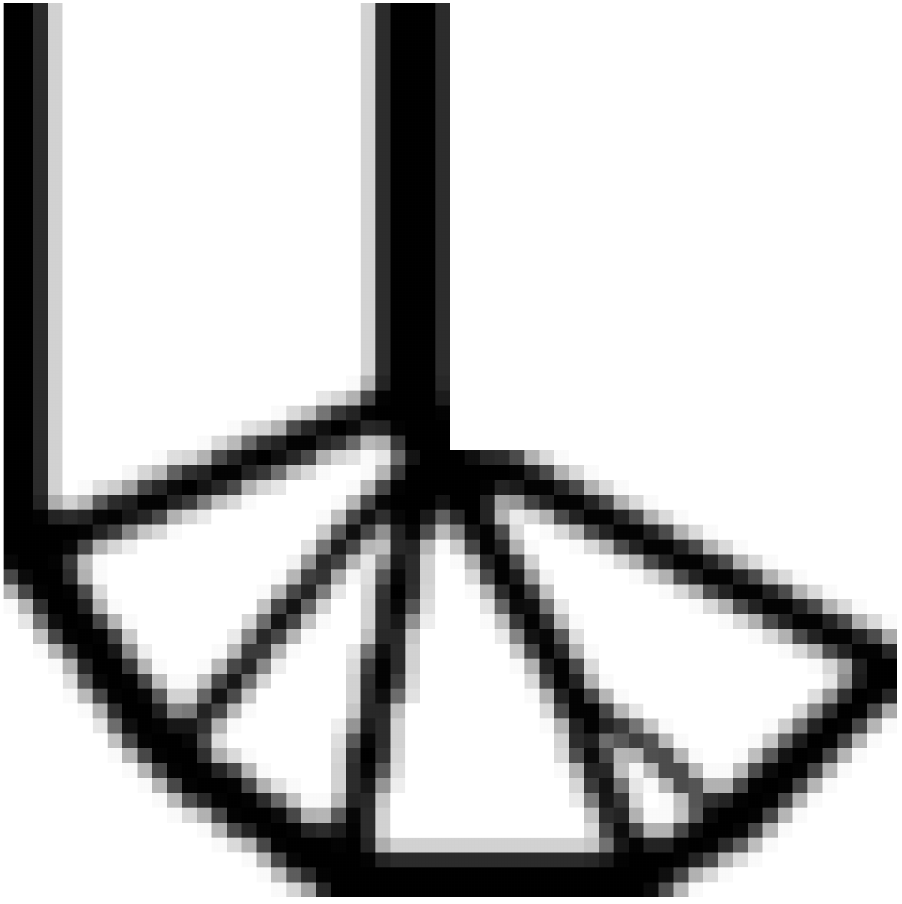} \\
        \bottomrule
    \end{tabular}
    \label{tbl:RBTO_L}
    \vspace{10px}
    \caption{The L-shaped beam: DTO results}
    \begin{tabular}{lcccc}
        \toprule
        $E$  & $1.15$ & $1.25$ & $1.35$ \\
        \midrule
        DTO 
        & \includegraphics[scale=0.43]{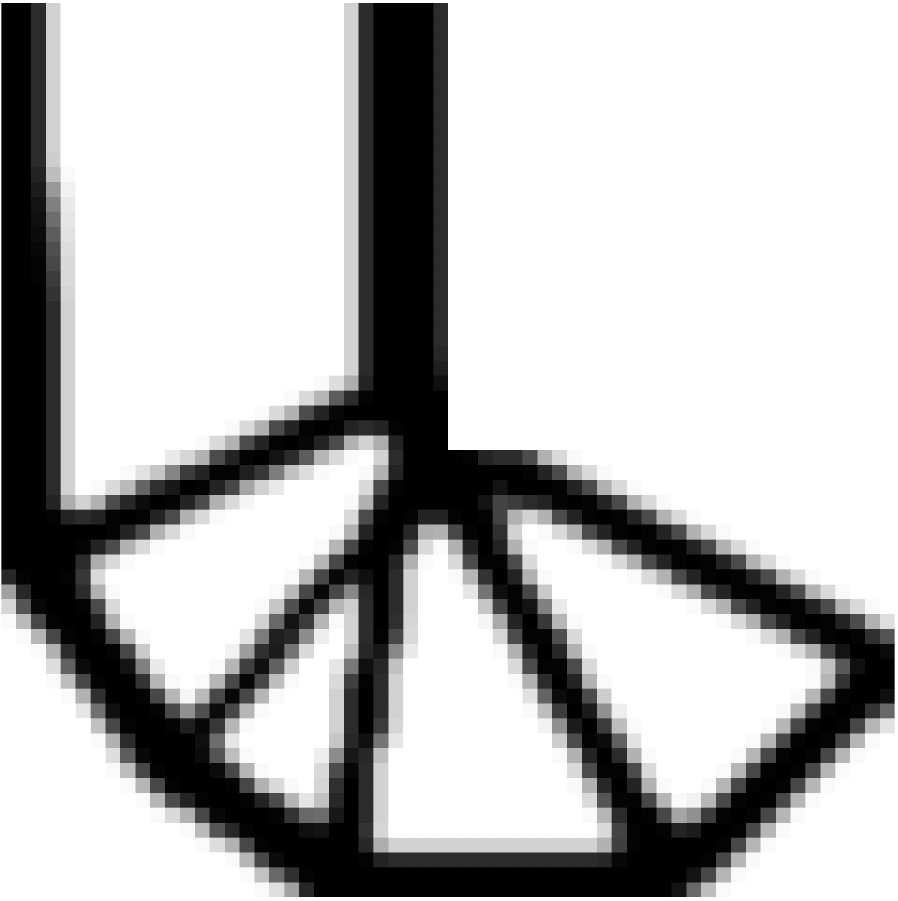}
        & \includegraphics[scale=0.43]{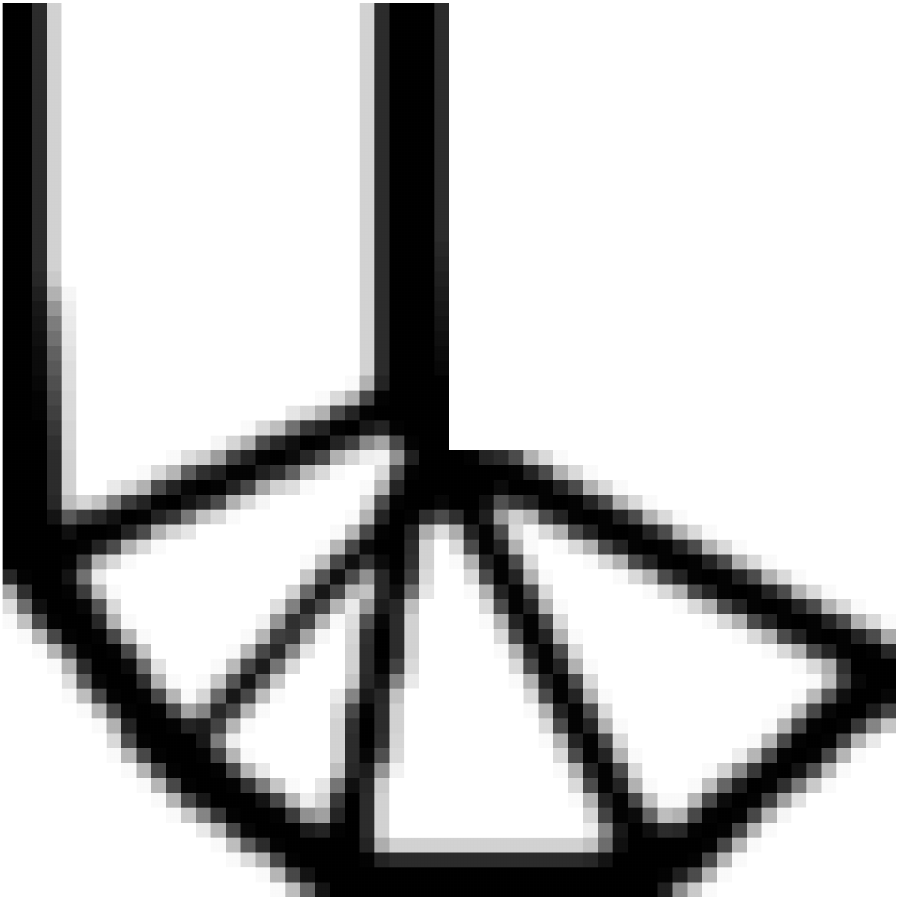}
        & \includegraphics[scale=0.43]{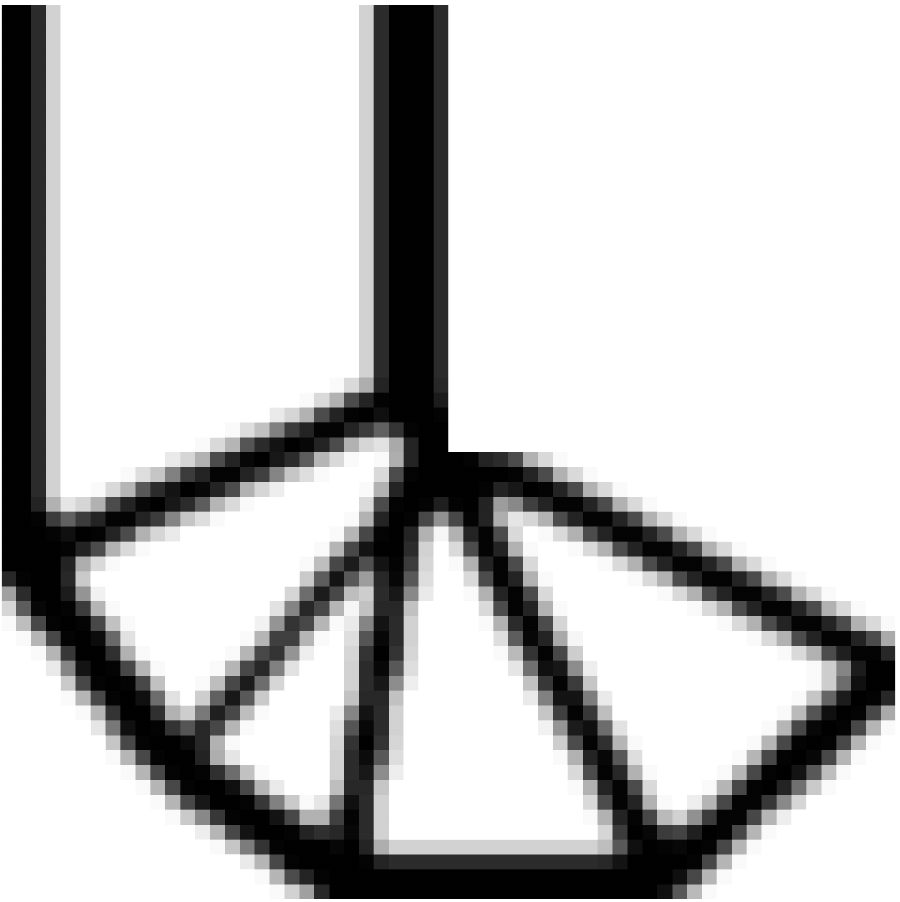} \\
        \midrule
        \multicolumn{1}{m{4em}}{Volume fraction} & 0.3542 & 0.3281 & 0.3097 \\
        \bottomrule
    \end{tabular}
    \label{tbl:DTO_L}
\end{table*}
\begin{table*}
    \centering
    \caption{The L-shaped beam: MCS, SRSM, and volume fraction}
    \begin{tabular}{lllllllll}
        \toprule
        & & & & \multicolumn{3}{c}{MCS} & \multicolumn{2}{c}{SRSM} \\
        \cmidrule(lrr){5-7}
        \cmidrule(lr){8-9}
        $\beta$ & Expected $P_f$ & $(a,b)$ & Volume fraction & $\mu$ & $\sigma$ & $P_f$ & $\mu$ & $\sigma$ \\
        \cmidrule{1-9}
        \multirow{3}{*}[-1em]{2} & \multirow{3}{*}[-1em]{0.02275} & $(1,1.3)$ & 0.3552 & 99.628 & 0.1844 & 0.02174 & 99.628 & 0.1844 \\
        \cmidrule{3-9}
        & & $(1,1.5)$ & 0.3291 & 99.420 & 0.2873 & 0.02198 & 99.420 & 0.2873 \\
        \cmidrule{3-9}
        & & $(1,1.7)$ & 0.3117 & 99.254 & 0.3688 & 0.02200 & 99.254 & 0.3688 \\
        \cmidrule{1-9}
        \multirow{3}{*}[-1em]{3} & \multirow{3}{*}[-1em]{0.001349} & $(1,1.3)$ & 0.3557 & 99.444 & 0.1832 & 0.001140 & 99.444 & 0.1832 \\
        \cmidrule{3-9}
        & & $(1,1.5)$ & 0.3295 & 99.127 & 0.2867 & 0.001120 & 99.127 & 0.2867 \\
        \cmidrule{3-9}
        & & $(1,1.7)$ & 0.3153 & 98.881 & 0.3664 & 0.001120 & 98.881 & 0.3664 \\
        \bottomrule
    \end{tabular}
    \label{tbl:L_Data}
\end{table*}
The second numerical example to illustrate our proposed algorithm is the L-shaped beam as shown in Fig.~\ref{fig:L}. The beam is fixed on its topmost edge and subject to a unit vertical load at the midpoint of the rightmost edge. As in the previous example, the volume of the beam is minimized while the probabilistic constraint is to limit the vertical displacement of the load application point (point B in Fig.~\ref{fig:L}). The maximum allowable vertical displacement of point B is chosen as $u^0=100$. Again, this number is only to make sure the results are visually friendly. The design domain is discretized using a $60\times60$ mesh of finite elements. One quarter of the mesh is removed to make the domain L-shaped by setting the element densities in this region to $0.001$ before proceeding to the next loop.

Compared to the MBB beam, this example uses a larger number of elements (3600 vs. 1200) resulting in much more time to finish the MCS. On the same machine, it took approximately 15 minutes to produce one row in Table~\ref{tbl:MBB_Data} but nearly 2.5 hours for the L-shaped beam. Furthermore, after examining the results, almost the same trends as in the previous example are observed. Thus, we decide to run MCS for only six cases as shown in Table~\ref{tbl:L_Data}. Those cases are the combinations of $\beta=\{2,3\}$ and $(a,b)=\{(1,1.3),(1,1.5),(1,1.7)\}$. The failure probabilities of the vertical displacement constraint at point B, and the statistical moments of that point's vertical displacement are calculated using both the MCS and SRSM in Table~\ref{tbl:L_Data}. Table~\ref{tbl:RBTO_L} and \ref{tbl:DTO_L} show the RBTO and DTO results of the L-shaped beam, respectively. Lastly, the CDF of point B vertical displacement and its last 10 points are displayed respectively in Fig.~\ref{fig:L_CDF_all} and \ref{fig:L_CDF_10}. Several insights and discussions are provided in the next section.

\section{Discussions}
\label{disu}
The results from the two examples are studied in this section for comparison, verification, and insights.

Visually, it is hard for human eyes to detect any differences between the RBTO results and the corresponding DTO results using the mean values of the Young's modulus. We are also not aware of any tools to find such differences, at least in the field of topology optimization, and hence believe this would be a research direction of great utility, especially when considering uncertainty and different material models. In the meantime, to compare those results qualitatively, readers can compile them column-wise (i.e., $(a,b)=(1,1.5)$ in Table~\ref{tbl:RBTO_MMB} and $E=1.25$ in Table~\ref{tbl:DTO_MMB}) into short animations, which are not possible to show in here. It turns out that there are material redistribution among most of the results, thickening or thinning of features, and even removal or addition of features (i.e., $(a,b)=(1,1.3)$ in Table~\ref{tbl:RBTO_MMB} vs. $E=1.15$ in Table~\ref{tbl:DTO_MMB}; $(a,b)=(1,1.7)$ in Table~\ref{tbl:RBTO_L} vs. $E=1.35$ in Table~\ref{tbl:DTO_L}).
\begin{figure}
    \centering
    \subfloat[ The last 10 points of the CDF]{\includegraphics[width=0.45\textwidth]{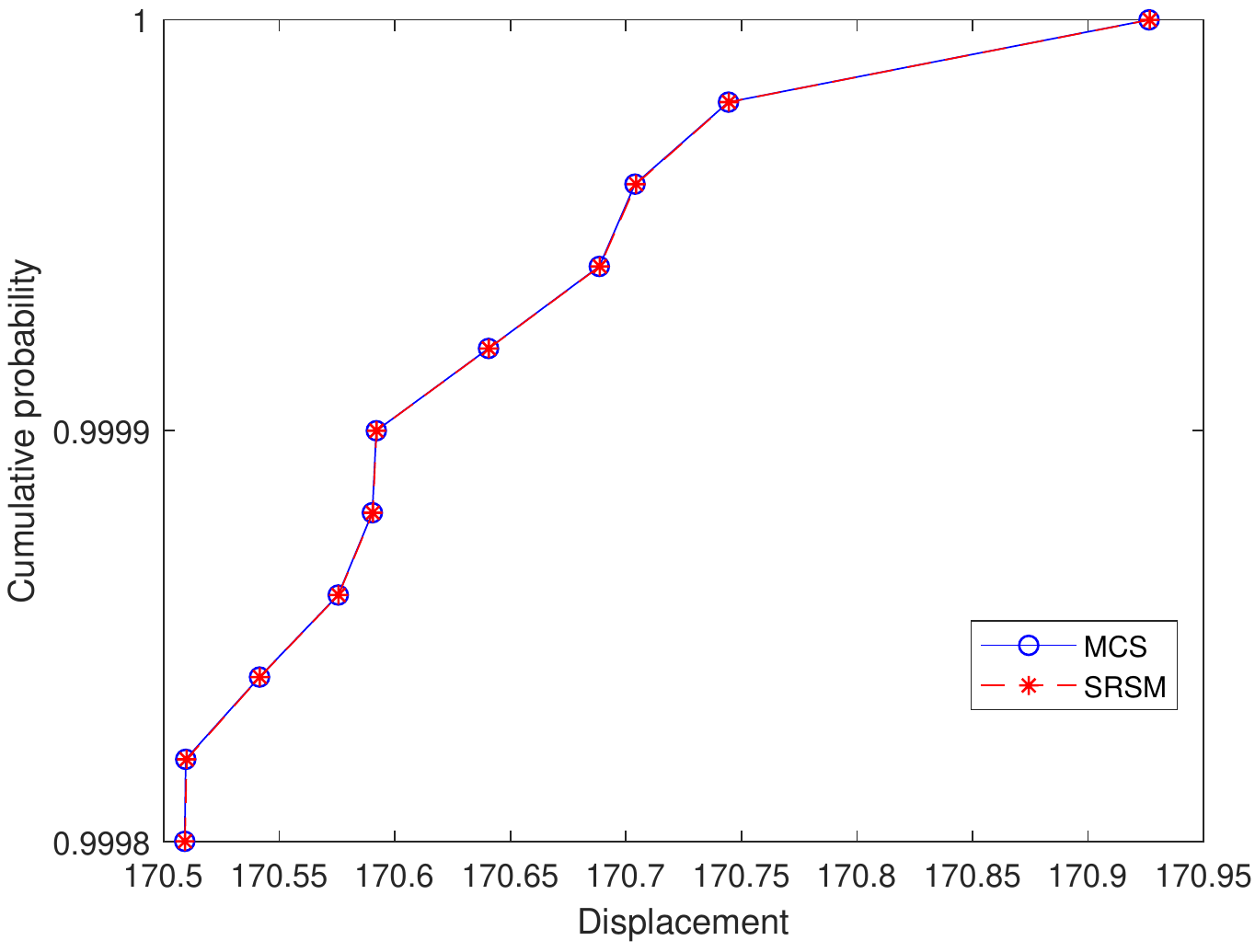}\label{fig:MBB_CDF_10}} \\
    \subfloat[The full CDF plot]{\includegraphics[width=0.45\textwidth]{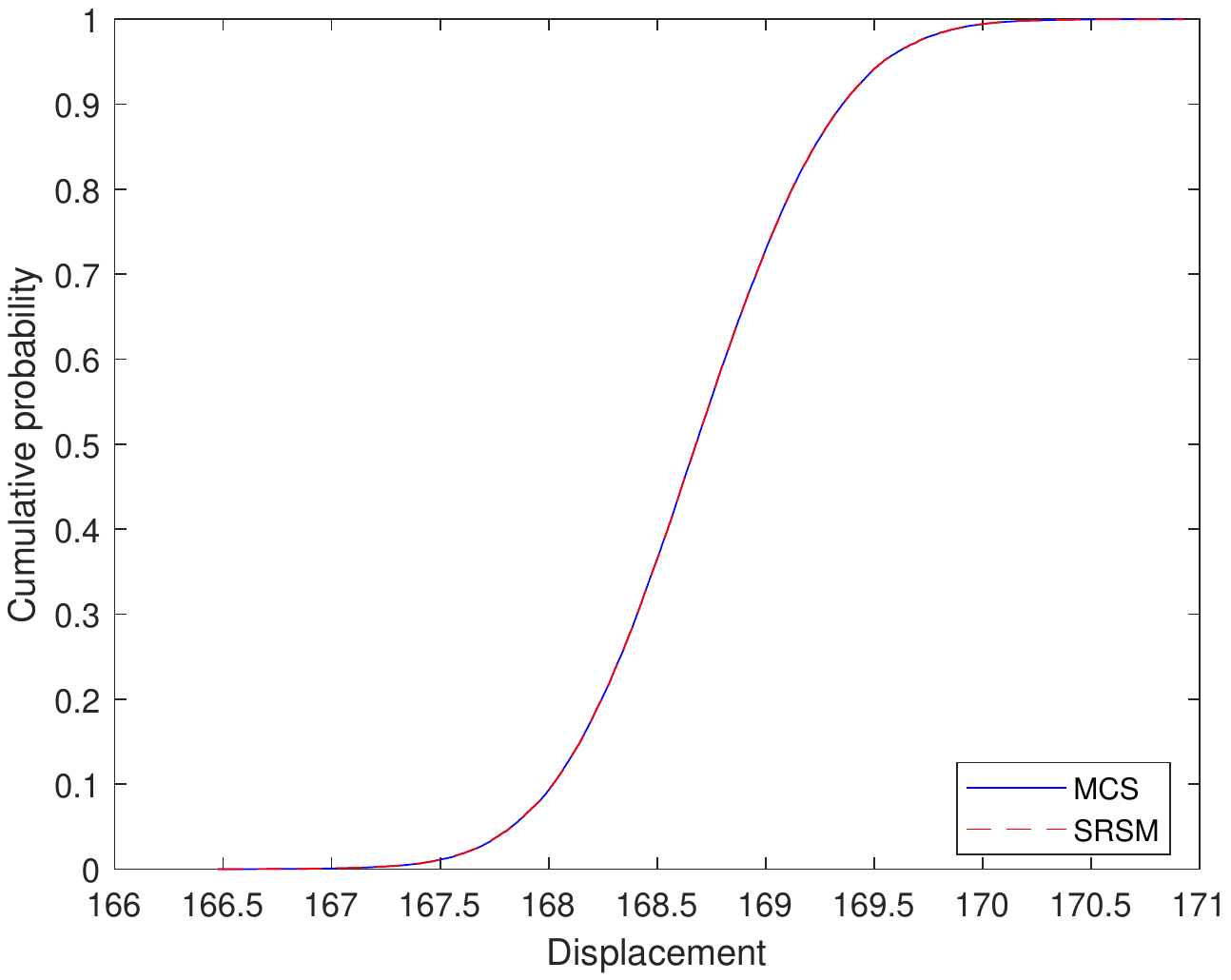}\label{fig:MBB_CDF_all}}
    \caption{The MBB beam: CDF plot for $\beta=2.5$ and $(a,b)=(1,1.5)$}
\end{figure}

The optimized designs together with the numbers in Table~\ref{tbl:MBB_Data} and \ref{tbl:L_Data} reveal several trends:
\begin{enumerate}
\item For the same reliability target, increasing variability of the material property (i.e., bigger $b$ in $(a,b)$) will decrease the volume fraction. This is reflected in thinner features or complete removal of features. For example, the features of the four designs in the first row of Table~\ref{tbl:RBTO_MMB} become thinner when increasing the range, making their ``inner'' spaces sparser. Smaller features can be translated into smaller cross section (i.e., of beam or truss member). This trend is also visible if ones look at the DTO results, which are found using the mean of the range $(a,b)$. If MCS are run on the DTO results, the failure probabilities will be around $0.5$, which is equivalent to $\beta=0$. This trend can be explained by the fact that the Young's modulus in the two examples is modelled by the uniform distribution, in which all intervals of the same length on the distribution's support have the same probability, and displacement is smaller with bigger value of the Young's modulus. Also, the mean and standard deviation of the vertical displacement have opposite trends: the mean decreases while the standard deviation becomes bigger. This is obviously caused by the increasing range while keeping the maximum allowable displacement $u^0$ constant.
\item For the same range, a higher reliability target will increase the volume fraction. The changes in most cases are pretty small, leading to hard-to-detect differences in the optimized designs. However, such small changes achieve significantly distinct failure probabilities if viewing the second column of Table~\ref{tbl:MBB_Data} and \ref{tbl:L_Data}. Possibly smaller failure probability enforces tighter bounds on the vertical displacement (i.e., both the mean and standard deviation becomes smaller), which in turn requires a volume increment.
\begin{figure}
    \centering
    \subfloat[The last 10 points of the CDF]{\includegraphics[width=0.45\textwidth]{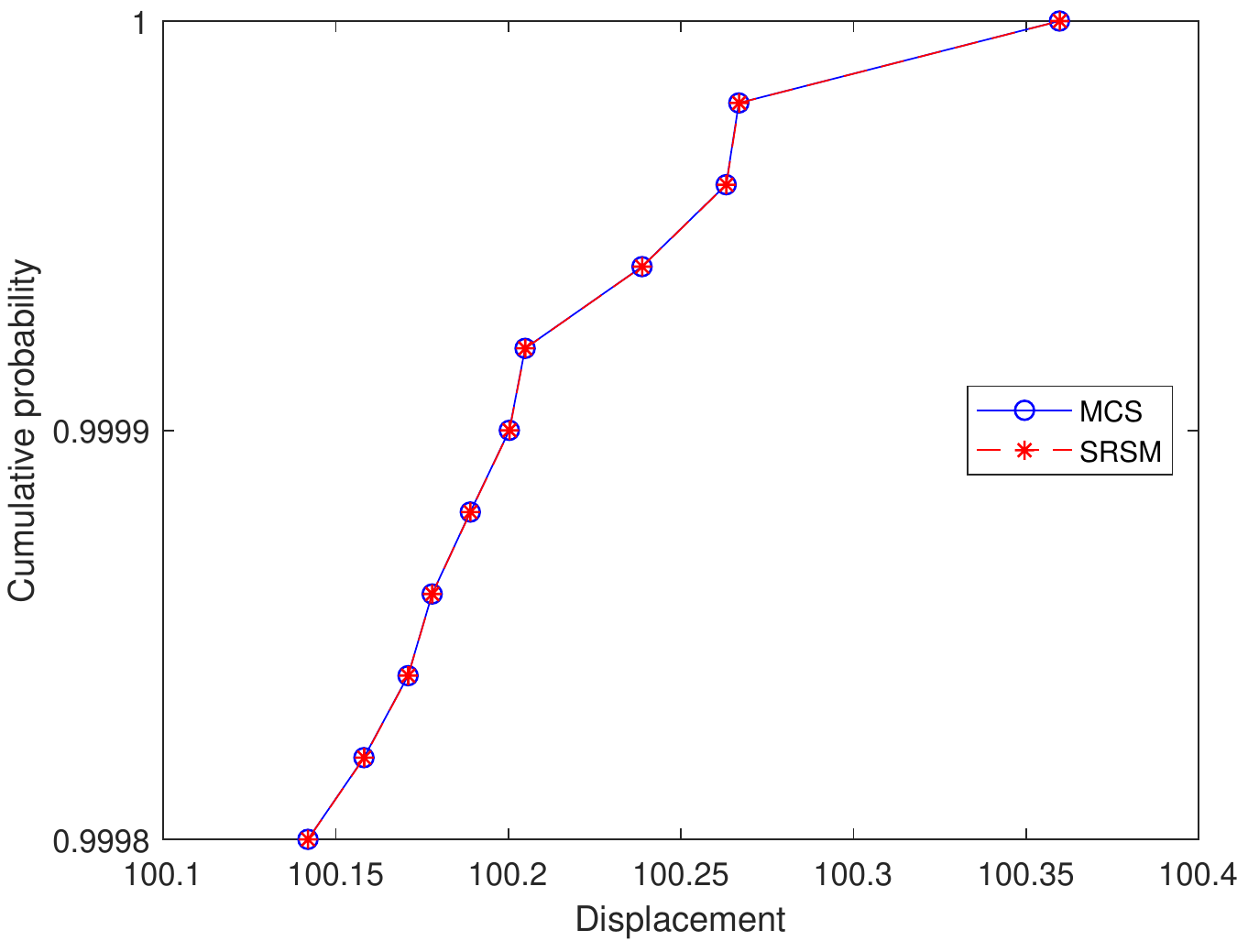}\label{fig:L_CDF_10}} \\
    \subfloat[The full CDF plot]{\includegraphics[width=0.45\textwidth]{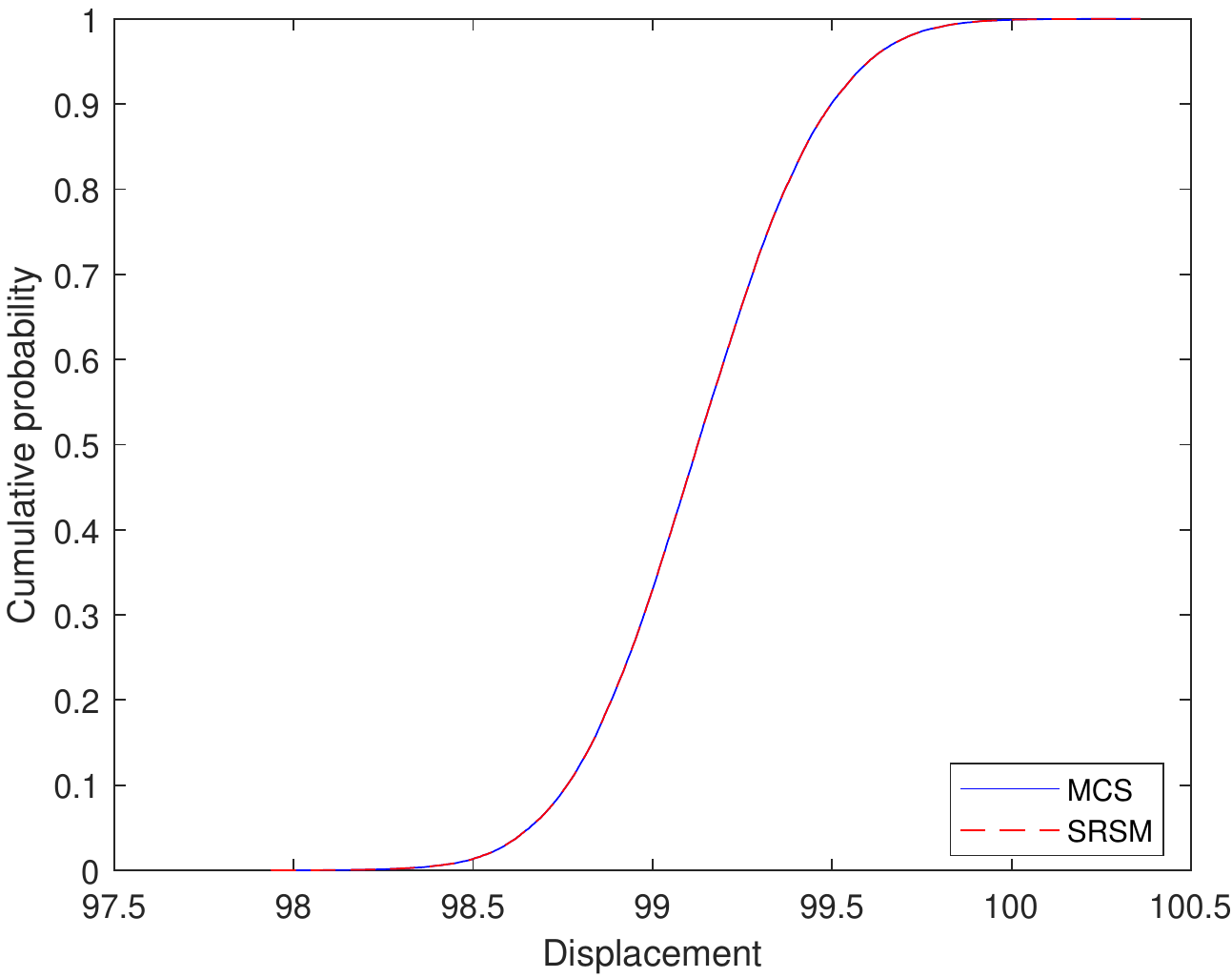}\label{fig:L_CDF_all}}
    \caption{The L-shaped beam: CDF plot for $\beta=3$ and $(a,b)=(1,1.5)$}
\end{figure}
\item The failure probabilities computed by the MCS (the $7^{\text{th}}$ columns in Table~\ref{tbl:MBB_Data} and \ref{tbl:L_Data}) in most cases are smaller than the corresponding expected values (the $2^{\text{nd}}$ columns), which proves the correctness of our proposed algorithm. However, there are three outliers in Table~\ref{tbl:MBB_Data} ($\beta=2$ and $(a,b)=(1,1.7)$; $\beta=3$ and $(a,b)=\{(1,1.5),(1,1.7)\}$). Those values are somewhat expected. First, there are many levels of approximation in our algorithm from the KL expansion to inverse reliability analysis and the SRSM, which accumulates error in the final results. Second, the random data points used in the MCS, which are generated by Latin hypercube sampling, are scattered uniformly over the whole distribution of the vertical displacement. The results would be much more accurate if there are more data points in the tail of the distribution, which usually has very small probability and where failure happens (e.g., using importance sampling). Third, the pseudorandom number generator in Matlab works by taking an initial seed, which is 0 for the results in Table~\ref{tbl:MBB_Data} and \ref{tbl:L_Data}, and then generating a deterministic sequence of numbers. Different seeds produce different sequences resulting in slightly different probabilities. It is quite likely that for the seed we chose, an ``unfavorable'' sequence is accompanied leading to those outliers. Based on the above arguments, we are confident that our algorithm is working correctly.
\item The MCS results in all cases agree well with the SRSM ones. From our observation, it is only possible to see the difference in the mean and standard deviation with 6 or more significant figures. Because of such small difference, the full CDF plots of the MCS and SRSM in Fig.~\ref{fig:MBB_CDF_all} and \ref{fig:L_CDF_all}, as well as the last 10 points of those CDFs in Fig.~\ref{fig:MBB_CDF_10} and \ref{fig:L_CDF_10} are almost indistinguishable. Thus, in the settings of out algorithm the SRSM works well and can be a strong, much cheaper alternative to the MCS. 
\end{enumerate}

\section{Conclusions}
\label{concl}
In this work, we have presented a comprehensive approach to reliability-based design optimization considering random field uncertainty, specifically for application to topology optimization. The approach follows the Sequential Optimization and Reliability Assessment (SORA)  approach formulated for random variable uncertainty, but extends the approach to consider random field uncertainty enabled using a Karhunen$-–$Lo\`{e}ve expansion and stochastic response surface. We have demonstrated the approach on two simple design examples and have shown that our framework is an efficient method to reliability-based design within a topology optimization approach to structural design. We have validated the approach by comparing the results from the SORA method to a Monte Carlo simulation.

The proposed approach is a first step to a comprehensive approach for design for additive manufacturing, or 3-D printing.  Additive manufacturing approaches can lead to uncertainty in material properties, such as the modulus of elasticity, since the material is generally deposited in a sequential fashion. This material uncertainty can be represented as random field uncertainty as we have presented in this paper. Future work towards realizing topology optimization as a comprehensive approach to design for additive manufacturing includes several additional research directions. In general, many additive manufacturing processes deposit material in a layered pattern, so consideration of material anisotropy, in addition to material uncertainty, should be explored as a key feature of topology optimization design approach. Additionally, many additive manufacturing applications use polymeric materials which do not follow linear elastic material assumptions, but rather are better modeled assuming hyperelastic or viscoelastic material assumptions. Inclusion of material non-linearly may influence the design resulting from the TO process. Also, while we have considered uncertainty in the optimization constraints, leading to reliability-based design optimization, one should also consider uncertainty in the optimization objective function, leading to robust reliability based design optimization. In terms of validation, application of the proposed approach to more realistic design cases should be conducted. The simple examples in this paper were chosen to best illustrate the proposed SORA approach to topology optimization; however, design cases such as the design of airless tires, or tweels, should be conducted to better understand how the introduction of material uncertainty influences the resulting design. 

\small{
\noindent\textbf{Replication of Results} \\
Matlab code for the two examples available at the Design Engineering Lab Github site \url{https://github.com/DesignEngrLab} in the project titled \textbf{RBTO}.

\noindent\textbf{Compliance with Ethical Standards} \\
The authors declare they have no conflict of interest.
}
\begin{acknowledgements}
The first author would like to thank the Vietnam Education Foundation (VEF) for their financial support through the VEF fellowship, and Prof. Krister Svanberg for his assistance on the MMA code.
\end{acknowledgements}

\bibliographystyle{spbasic}      %
\bibliography{refs}   %

\end{document}